# BULK DIFFUSION IN A SYSTEM WITH SITE DISORDER

By Jeremy Quastel[1]

*University of Toronto*

We consider a system of random walks in a random environment interacting via exclusion. The model is reversible with respect to a family of disordered Bernoulli measures. Assuming some weak mixing conditions, it is shown that, under diffusive scaling, the system has a deterministic hydrodynamic limit which holds for almost every realization of the environment. The limit is a nonlinear diffusion equation with diffusion coefficient given by a variational formula. The model is nongradient and the method used is the "long jump" variation of the standard nongradient method, which is a type of renormalization. The proof is valid in all dimensions.

**1. Introduction.** Consider a system of particles occupying sites of a multidimensional integer lattice. The particles are attempting jumps to nearest neighbor sites at rates which depend on both their position and the objective site. The rates themselves come from a quenched random field, and are chosen so that the system satisfies a detailed balance condition with respect to a family of random Bernoulli measures (the random field Ising model at infinite temperature). The interaction between the particles is given by a hard core exclusion rule: Attempted jumps to occupied sites are simply suppressed.

Such systems have been used to model electron transport in doped crystals [1, 2, 3, 8, 10, 12, 14]. In this case, the hard core exclusion rule is given by the Pauli exclusion principle. The crystal itself creates a periodic field in which the particles move. In the presence of impurities, the field is random.

The purpose of this article is to study the transport properties of such a system. In particular we are interested in the influence of the random field on the rate of bulk diffusion. In the absence of the exclusion rule, this is simply the result of diffusive scaling of a single particle moving with

Received April 2004; revised July 2005.
[1]Supported in part by NSERC.

reversible rates in the random field. At the other extreme, if the exclusion is present but the field is constant, the bulk diffusion also turns out to be independent of the density of particles. However, when both the random field and the exclusion rule are present, one sees a nonlinear dependence of the bulk diffusion on the density. This diffusion coefficient can be computed by the Green–Kubo formula, but we have not been able to find such a computation in the literature. We give a variational formula for the bulk diffusion, which is equivalent to the Green–Kubo formula, and which we can establish rigorously.

We should point out that in principle it is not even clear that under diffusive scaling such a system has a hydrodynamic limit. The main work here is to show that in fact all the influences of the random field are, on a sufficiently large scale, contained in a diffusion coefficient, which depends on the statistics of the field, but not on the randomness itself. One of the main steps in this direction is to show that the system in a box of side length $L$ has a spectral gap no smaller than $cL^{-2}$.

To establish the hydrodynamic limit, one needs to prove some version of Fick's law, namely, to replace the microscopic current by the gradient of the density field multiplied by the diffusion coefficient. The system we are considering turns out to be of nongradient type. Roughly speaking, the gradient condition says that the microscopic current is already of gradient form. At this time, the only method for the nongradient systems is the one developed by Varadhan. The idea is to replace the current by a gradient plus a fluctuation term. However, in our case, such a decomposition cannot hold microscopically, because the fluctuations of the gradient of the density field arising from the random field are very large, and one can only make sense of these gradients over sufficiently large mesoscopic distances where stochastic fluctuations are reduced by the central limit theorem. It is necessary to perform the decomposition in some mesoscopic scale.

The model we are considering has been studied in the physics literature by means of rough approximations and Monte Carlo simulations. The case of a period two field in one dimension has been solved exactly [20, 24]. Surprisingly, the resulting equation is linear. This special case turns out to be an "almost gradient system," by which we mean that the fluctuation term is explicit. In other words, one can find the minimizer in the variational problem for the diffusion coefficient. This is certainly not the case for other fields, and it can be shown that the diffusion coefficient is nonconstant in general.

The main result of this paper was described in detail in [16] along with a sketch of the proof and an unpublished, incomplete manuscript [18]. An important motivation for the method described there is that the nongradient method as developed in [17, 22, 23] does not work in low dimensions. In higher dimensions, such an approach can be made to work by subtracting a term from the microscopic gradient term to make it mean zero. One then has the nontrivial problem to show that the subtracted term vanishes in the limit, and this can be done if $d \geq 3$. This was suggested in [16], and a proof following these lines appeared recently [6]. This inspired us to write up the details of our unpublished manuscript carefully as the present paper. The advantage of the original approach in [16] is that it works in all dimensions and elucidates the connection to renormalization. The present article is independent of [6]; the only result we have used beyond [16] and [18] is the proof of the spectral gap of the Bernoulli–Laplace version of the model by [4], which improves on what was written in [18]. Note that both [4] and [6] need the moving particle lemma (Lemmas 5.1 and 5.2 of this article), which is in some sense what tells you that the system is diffusive.

In addition, we prove here the continuity of the diffusion coefficient in the full interval $[0,1]$ and the uniqueness of the hydrodynamic equation. The uniqueness needed to be added as an assumption in [6], and the continuity was proven only on $(0, 1)$. Note that the continuity of the diffusion coefficient at 0 is of considerable interest because there it becomes the asymptotic diffusivity of a single particle in the medium, a classical homogenization problem with a different variational formula that has to be related to ours.

The review article [16] contains a sketch of the main idea of the proof contained in this article. The sketch contains a small misprint which has unfortunately led to a great deal of confusion. The formula for a certain central limit theorem variance,

$$(1.1) \quad \limsup_{K \to \infty} \mathbb{E}[\langle K^{-1} W_K, (-K^2 L_K)^{-1} K^{-1} W_K \rangle] = C m^2 (1 - m^2) \sigma^{-1},$$

is stated without a numerical factor $C$ on the right-hand side. We would like to explain the nature of this minor mistake, which does not affect the integrity of the proof. In fact, as long as it is finite, the exact value of the constant $C$ is not even relevant in the proof. To explain this, we need to describe the basic idea of the proof. This is a bit easier in one space dimension.

The main problem in proving the hydrodynamic limit is to understand the large scale behavior of the average current $w_{x,x+1}$ between site $x$ and $x+1$, multiplied by a factor $\varepsilon^{-1}$, the size of the system. The current is a function of the particle numbers and field values at those two sites. One is allowed to take some average of these currents over a larger, but still microscopic block. Write the sum on the block as $\bar{w}_{x,y} = \sum_{z=x}^{y-1} w_{z,z+1}$. Averaging does not appear to help too much unless $w_{x,x+1}$ were of *gradient* form, $\tau_{x+1} h - \tau_x h$ for some $h$. Then there would be an easy term by term cancellation. But in this model the current is not of gradient form unless $\alpha = 0$. It might work, however, if we could find an $h$ such that the difference can be made appropriately small. This is the standard approach to nongradient systems, which as explained, works for this model only in higher dimensions. Our alternative approach is to use a *long jump* version of the current, $w_{x,y}$, which is the same function evaluated with variables at sites $x$ and $y$, but now $x$ and $y$ are chosen far apart. It turns out to be true that one can replace $\bar{w}_{x,y}$ by $Cw_{x,y}$, where $C$ is related to the diffusion coefficient. This helps somewhat because one needs to normalize the sum in $\bar{w}_{x,y}$, and a normalization factor of $y - x$ in the denominator is now available. But to compensate for the huge $\varepsilon^{-1}$ factor in front of the current, one would have to take $y - x = \varepsilon^{-1}$ which is no longer microscopic. So one needs something more. One can get the needed extra factor by averaging the long jump ends $x$ and $y$ over some large microscopic blocks, say, $x \in \Lambda_1$ and $y \in \Lambda_2$, where $\Lambda_1$ and $\Lambda_2$ are two blocks of side length $K$ beside each other. Then one gets a free factor $K^{-1/2}$ ($K^{-d/2}$ in dimension $d$) from the central limit theorem, and together with the $K^{-1}$ from the typical distance between $x$ and $y$, this gives a prefactor $K^{-(d+2)/2}$ which can beat the $\varepsilon^{-1}$ if $K$ is chosen as large as $\varepsilon^{-2/(d+2)}$.

So we want to replace

$$(1.2) \quad K^{-1} Av_{\substack{x \in \Lambda_1 \\ y \in \Lambda_2}} \sum_{z=x}^{y-1} w_{z,z+1} \quad \text{by } K^{-1}W_K = K^{-1} Av_{\substack{x \in \Lambda_1 \\ y \in \Lambda_2}} w_{x,y}.$$

If one instead tries to replace $Av_{z,z+1 \in \Lambda_1 \cup \Lambda_2} w_{z,z+1}$ by $K^{-1}W_K = K^{-1} \times Av_{\substack{x \in \Lambda_1 \\ y \in \Lambda_2}} w_{x,y}$, the resulting central limit theorem variance will not vanish. This is not a surprise and is simply because one forgot the extra average on the left-hand side of (1.2), so the corresponding variables are weighted differently than those on the right-hand side of (1.2). But it appears some readers drew the erroneous conclusion that this is because of the missing factor $C$ in (1.1). We regret the misprint and the confusion it has caused, but it does not disqualify in any way the main idea of the proof as described in [16] and [18].

The proof of the hydrodynamic limit contained in this article is essentially the same as that of [16] and [18]. The main new material is the proof of the continuity of the diffusion coefficient, the uniqueness of the hydrodynamic equation and the proof of (1.1), the details of which were not included in either [16] or [18]. This is in Sections 7 and 8. It requires a nontrivial extension of the standard nongradient computations, as one is missing in this problem the usual average over translations.

The paper is organized as follows. In Section 2 we describe the model and the main results. In Section 3 we give a sketch of the proof of the hydrodynamic limit, emphasizing the new problems that arise because of the random field. Sections 4–8 contain the details of the proof. In Section 4 we recall some standard results in perturbation theory which allow us to use the variance method for nongradient systems on functions whose range is up to a small constant times $\varepsilon^{-2/(d+2)}$ as long as they satisfy an "integration by parts" estimate. In Section 4 we prove the key "moving particles lemma" which gives the spectral gap and the two-block estimate and the "integration by parts" estimate for the long jump current. Sections 6–10 are devoted to the computation of asymptotic variances which is the heart of the nongradient method. The standard material is in Sections 6 and 7. In Sections 8 and 9 the computation is extended to the "long jump current" where one no longer has an average over shifts. In Section 10 we complete the proof by showing that the long jump current can be replaced by its average with respect to local equilibrium at the scale $\varepsilon^{-2/(d+2)}$ where the fluctuations have been sufficiently dampened. Section 11 contains the results about continuity and a type of Hölder continuity of the diffusion coefficient and Section 12 the uniqueness of the hydrodynamic equation.

**2. The model.** Let $\alpha_x$, $x \in \mathbb{Z}^d$, be a bounded ($|\alpha_x| \leq B$) stationary, ergodic random field satisfying the following mixing conditions: For some $\gamma \geq \max\{4, 2(d+2)/d^2\}$, there is a constant $C < \infty$ such that, for all $l > 0$, for all $f \in \mathcal{F}_{\Lambda_l}$ with $\mathbb{E}[f] = 0$,

$$(2.1) \qquad \mathbb{E}[|Av_{x \in \Lambda_K} \tau_x f|^\gamma] \leq C(K^{-1}l)^{\gamma d/2} \mathbb{E}[f^\gamma],$$

where $\tau_x$ is the shift by $x$, $\Lambda_K$ is the box of side length $K$, and $\mathbb{E}$ is expectation with respect to the field. Typical examples are $\alpha$ periodic, or $\alpha_x$ could be independent, identically distributed random variables. But (2.1) is quite general. For example, one could divide $\mathbb{Z}^d$ into boxes of side length $A$. On each of the boxes one could pick at random from a finite list of patterns (functions on $\{1, \ldots, A\}^d$).

For each $\varepsilon = L^{-1}$, $L$ a positive integer, we have a system of $N = O(\varepsilon^{-d})$ particles on $\varepsilon \mathbb{Z}^d / \mathbb{Z}^d$ moving in this field. At most one particle is allowed at each site. A particle at $x$ attempts to jump to nearest neighbor sites $y$ at rate

$$(2.2) \qquad \varepsilon^{-2}(1 + e^{\alpha_y - \alpha_x}).$$

If there is no particle in the way, the particle is allowed to jump. However, if there is a particle in the way, the jump is suppressed, and everything starts again. All the particles are doing this independently of each other, and since time is continuous, one can ignore the occasion of two particles trying to jump onto each other simultaneously.

The state space of our process is $\{0,1\}^{\varepsilon \mathbb{Z}^d / \mathbb{Z}^d}$. Configurations are denoted $\eta$. $\eta_x = 1$ or $0$ depending on whether there is or is not a particle at $x$. Our system is a Markov process on this state space with a generator given by $\varepsilon^{-2} L_\varepsilon$, where

$$(2.3) \qquad L_\varepsilon = \sum_{|x-y|=1} a_{xy}(\eta) \nabla_{xy},$$

where the rates are given by

$$a_{xy}(\eta) = 1 + e^{-(\alpha_y - \alpha_x)(\eta_y - \eta_x)}, \tag{2.4}$$

and the lattice gradient is given by

$$\nabla_{xy} f(\eta) = f(T_{xy}\eta) - f(\eta), \tag{2.5}$$

where $T_{xy}\eta$ represents $\eta$ with the occupation numbers at $x$ and $y$ exchanged. For each $\lambda$, the generator is reversible (self-adjoint) with respect to the product measure

$$Z^{-1} \exp \sum_{x \in \varepsilon \mathbb{Z}^d/\mathbb{Z}^d} (\alpha_x + \lambda)\eta_x. \tag{2.6}$$

$Z$ is the normalization. The parameter $\lambda$, which is called the chemical potential, can be adjusted to vary the average density of particles. The relation between $\lambda$ and the density, $m$, is as follows: For $0 \leq m \leq 1$, $\lambda$ is chosen so that

$$\mathbb{E}\left[\frac{e^{\alpha_0 + \lambda}}{1 + e^{\alpha_0 + \lambda}}\right] = m. \tag{2.7}$$

We could alternatively fix the number of particles so that the density was $m$. Then the process is reversible and ergodic with respect to that measure

$$Z^{-1} \exp \sum \alpha_x \eta_x \; \delta_{\{Av_x\eta_x = m\}}. \tag{2.8}$$

Taken over the allowable $m = i\varepsilon^d$, $i = 0, \ldots, \varepsilon^{-d}$, these give us a full set of ergodic invariant measures. Note that the measures (2.8) are not simple to describe. The Dirichlet form is given by

$$D(f) = E[f, (-L)^{-1} f] = \sum_{x \sim y} E[(\nabla_{xy} f)^2]. \tag{2.9}$$

The empirical density field $\mu_\varepsilon \in M(\mathbb{T}^d)$, the set of measures on the $d$-dimensional torus $\mathbb{T}^d$, is given by

$$\mu_\varepsilon(d\theta) = \varepsilon^d \sum_x \eta_x \delta_{\varepsilon x}, \tag{2.10}$$

where $\delta_\theta$ gives mass one to the point $\theta \in \mathbb{T}^d$. Let us choose initial distributions of our process so that, for some given smooth $m^0 : \mathbb{T}^d \to [0,1]$, $\mu_\varepsilon \Rightarrow m^0(\theta) \, d\theta$ in probability. The corresponding probability measure on $D([0,T] \to M(\mathbb{T}^d))$, the space of right continuous trajectories with left limits in $M(\mathbb{T}^d)$, will be denoted $P_\varepsilon$. Of course, $P_\varepsilon$ depend on the field, $\alpha$.

We can now state our main result.

THEOREM 1. *For almost every realization $\alpha$ of the random field, $P_\varepsilon \Rightarrow \delta_{m(t,\theta)d\theta}$, as $\varepsilon \to 0$, where $m$ is the unique weak solution of*

$$(2.11) \qquad \frac{\partial m}{\partial t} = \nabla \cdot D(m)\nabla m, \qquad m(0,\theta) = m^0(\theta).$$

*The diffusion matrix $D(m)$ is a nonrandom continuous function on $[0,1]$. It is given by the following formulae. For $m \in (0,1)$,*

$$D(m) = \sigma(m)\lambda'(m),$$

*where the conductivity, $\sigma$, is a symmetric matrix whose associated quadratic form is given by the variational formula,*

$$(2.12) \quad (\beta, \sigma(m)\beta) = 2\inf_g \mathbb{E}\left[\left\langle \sum_{e=1}^d \left(\beta_e(\eta_e - \eta_0) - \nabla_{0,e}\sum_x \tau_x g\right)^2 \right\rangle_m\right].$$

*The infimum is taken over all local functions $g(\eta,\alpha)$ of the configuration and the field. The shift is given by $\tau_x g(\eta,\alpha) = g(\tau_x\eta,\tau_x\alpha)$. The expectation $\mathbb{E}$ is over the random field $\alpha$ and the expectation $\langle \cdot \rangle_m$ is over the infinite product measure (2.6), where $\lambda$ is chosen as in (2.7). For $m = 0$, $D(0)$ is the limiting covariance of a free particle in the field, given by the classical homogenization formula,*

$$(2.13) \quad (\beta, D(0)\beta) = 2z^{-1}\inf_{U(\alpha)} \mathbb{E}\left[\sum_e (e^{\alpha_0} + e^{\alpha_e})(\beta_e + \tau^{-e}U - U)^2\right],$$

*where $z = \mathbb{E}[e^{\alpha_0}]$. For $m = 1$, $D(1)$ is the limiting covariance of a test "hole" which coincides with the limiting covariance of a free particle in the field $\tilde{\alpha} = -\alpha$.*

The result was obtained earlier in [6] in $d \geq 3$ and under the assumption that the diffusion coefficient is continuous and that the limiting equation has a unique solution.

A few comments follow.

*Regularity of $D(m)$.* From (2.12), $\sigma(m)$ is upper semicontinuous for $m \in (0,1)$. It is not hard to check that $\lambda'(m)$ is continuous on $(0,1)$ and, therefore, $D(m)$ is upper semicontinuous on $(0,1)$. The test function $g \equiv 0$ in (2.12) shows that $\sigma(m) \leq Cm(1-m)$ for some $C < \infty$ and it is elementary to check that $\lambda'(m) \leq C/m(1-m)$ for another finite $C$. Hence, $D(m) \leq cI$ for some $c < \infty$ for all $m \in (0,1)$. From the moving particles Lemma 5.2, one obtains also a lower bound $D(m) \geq c^{-1}I$ for $m \in (0,1)$. In Section 9 we show that $D(m)$ is Hölder $1/2$ with a coefficient which may behave badly at the edge; $|D(m_1) - D(m_2)|^2 \leq C(m_1(1-m_1))^{-1}|m_1 - m_2|$ for some $C < \infty$. This follows from the characterization of the diffusion coefficient in Sections 6 and 7 without too much work. It is also shown there that $D(m)$ is continuous on the whole interval $[0,1]$.

*Weak solutions.* By a weak solution of (2.11), we mean a function $m : [0, T] \times \mathbb{T}^d \to [0, 1]$ satisfying

$$
(2.14) \quad \int_0^T \int_{\mathbb{T}^d} \frac{|\nabla m|^2}{m(1-m)} \, d\theta \, dt = \sup_\phi \left\{ \int_0^T \int_{\mathbb{T}^d} [m \nabla \cdot \phi - |\phi|^2 m(1-m)] \, d\theta \, dt \right\} < \infty
$$

for each $T < \infty$ and for smooth test functions $\phi$ on $\mathbb{T}^d \times [0, T]$,

$$
(2.15) \quad \int_{\mathbb{T}^d} \phi m \, d\theta \Big|_{t=0}^{t=T} = \int_0^T \int_{\mathbb{T}^d} \frac{\partial \phi}{\partial t} m \, d\theta \, dt - \int_0^T \int_{\mathbb{T}^d} \nabla \phi \cdot D(m) \nabla m \, d\theta \, dt.
$$

In Section 10 we show that, under the Hölder continuity proved in Section 9, such weak solutions are unique.

*Associated dynamics.* One can produce other dynamics with (2.8) as ergodic reversible measures. Consider the dynamics associated to the Dirichlet form

$$
(2.16) \quad \sum_{x,y} p_{y-x} E[b_{xy}(\eta)(\nabla_{xy} f)^2],
$$

where $p$ is finite range and symmetric and $b_{x,x+y}(\eta) = \tau_x b_{0,y}(\eta)$ for some finite range $b_{0,y}(\eta)$ bounded above and below (for $y$ in the range of $p$). The expectation is with respect to any one of the measures (2.8). The corresponding dynamics have a particle at $x$ attempting to jump to $x + y$ at rate $p_y(b_{0y}(\tau_x \eta) + b_{0y}(T_{0y} \tau_x \eta) e^{\alpha_{x+y} - \alpha_x})$. Since our methods are based on estimates involving the Dirichlet form, they extend easily to equivalent Dirichlet forms, and one obtains an identical theorem with (2.12) replaced by the infimum over

$$
(2.17) \quad 2 \mathbb{E} \left[ \sum_y p_y \left\langle b_{0y} \left( \beta \cdot y(\eta_y - \eta_0) - \nabla_{0,y} \sum_x \tau_x g \right)^2 \right\rangle_m \right].
$$

A natural example is $b_{0e} = e^{-\alpha_0} \eta_0 (1 - \eta_e) + e^{-\alpha_e} \eta_e (1 - \eta_0)$. The resulting dynamics have a particle at $x$ attempting to jump to each nearest neighbor site at rate $e^{-\alpha_x}$.

*Mixing conditions.* The mixing conditions given in (2.1) are chosen for convenience and are not meant to be optimal. They are mainly to point out that the theorem holds under some very weak mixing conditions on the variables $\alpha_x$. A nice problem is to consider the case of unbounded field $\alpha$ (the condition that $|\alpha_x| \leq B$ is very important for the method).

**3. Hydrodynamic limit.** In this section we give a sketch of the proof of the hydrodynamic limit. Many of the arguments are now standard and can be found, for example, in [11]. However, at some points new ideas are needed, especially in low dimensions where there is not sufficient averaging to control the fluctuations from the random field. We will sketch the approach, emphasizing where new methods are needed, leaving the rigorous proofs for Sections 4–8.

The evolution of the empirical density field $\mu_\varepsilon$ is described by the following set of stochastic integral equations:

$$\int_{T^d} \phi(t,\theta)\mu_\varepsilon(t,d\theta)\Big|_{t=0}^{t=T} = \int_0^T \int_{T^d} \frac{\partial \phi}{\partial t}(t,\theta)\mu_\varepsilon(t,d\theta)\, dt$$
$$- \int_0^T \varepsilon^{d-2} \sum_{|x-y|=1} (\phi(t,\varepsilon y) - \phi(t,\varepsilon x))w_{xy}\, dt$$
$$+ \int_0^T \varepsilon^d \sum_{|x-y|=1} (\phi(t,\varepsilon y) - \phi(t,\varepsilon x))\nabla_{xy}\eta\, dM_{xy}.$$

Here $M_{xy}$ are independent Poisson "sawtooth" martingales running at rates $\varepsilon^{-2}a_{xy}$ and the current

$$(3.1) \quad w_{xy} = a_{xy}(\eta_y - \eta_x) = \eta_y(1-\eta_x)(1+e^{\alpha_x-\alpha_y}) - \eta_x(1-\eta_y)(1+e^{\alpha_y-\alpha_x}).$$

Using the formula $E[|M_f(T)|^2] = \int_0^T E[Lf^2 - 2fLf]\, ds$ if $M_f(T) = f(\eta(T)) - f(\eta(0)) - \int_0^T Lf(\eta(s))\, ds$, true for any Markov process, one easily computes the quadratic variation of the martingale term,

$$(3.2) \quad E^{P_\varepsilon}\left[\left(\int_0^T \varepsilon^d \sum_{|x-y|=1} (\phi(t,\varepsilon y) - \phi(t,\varepsilon x))\nabla_{xy}\eta\, dM_{xy}\right)^2\right] \leq C\varepsilon^d,$$

where $C$ depends only on $T$ and $\varepsilon^d \sum_{x \in \varepsilon\mathbb{Z}^d/\mathbb{Z}^d} \phi^2(x)$.

If one starts with a nondegenerate product invariant measure ($E[\eta_0] \neq 0$ or 1), we have the equilibrium, or stationary process which we denote by $Q_\varepsilon$ and from the bound on the entropy of any initial distribution with respect to that reference measure, we obtain directly the bounds

$$(3.3) \quad H(P_\varepsilon/Q_\varepsilon) \leq C\varepsilon^{-d}, \qquad \int_0^T D(\sqrt{f_t})\, dt \leq C\varepsilon^{2-d},$$

where $f_t$ is the marginal density of the nonequilibrium process at time $t$. If $V$ is any bounded function, we can estimate by the Feynman–Kac formula, the spectral theorem for self-adjoint operators, and the variational formula for the principle eigenvalue of $L+V$,

$$T^{-1}\log E^{Q_\varepsilon}\left[\exp\left\{\int_0^T V(\eta(s))\, ds\right\}\right]$$

(3.4)
$$\leq \sup\left\{\langle Vf\rangle - \varepsilon^{-2} \sum_{|x-y|=1} \langle (\nabla_{xy}\sqrt{f}\,)^2\rangle\right\},$$

where the expectations are with respect to the invariant measure and the supremum is over densities relative to that measure. The superexponential bound

(3.5)
$$E^{Q_\varepsilon}\left[\exp\varepsilon^{-d}\sum_{x\in\mathbb{Z}^d/\varepsilon^{-1}\mathbb{Z}^d}\phi(\varepsilon x)\int_s^t w_{xx+e}(s)\,ds\right]$$
$$\leq 2\exp(t-s)\sum_{x\in\varepsilon\mathbb{Z}^d/\mathbb{Z}^d}\phi^2(x)$$

is obtained this way using the integration by parts formula

(3.6) $$\langle w_{xy}f\rangle = -\tfrac{1}{2}\langle(\eta_y-\eta_x)\nabla_{xy}f\rangle,$$

true for any $f$ and any invariant measure on any set containing $x$ and $y$. By Schwarz's inequality, we have

(3.7) $$\langle w_{x,x+e}f\rangle \leq C\langle(\nabla_{xy}\sqrt{f}\,)^2\rangle^{1/2},$$

which gives (3.5). Once one has (3.2) and (3.5), tightness of the measures $P_\varepsilon$ follows by a standard argument from Garsia's lemma (see [11] for details).

Once one has tightness, it remains to identify the limit measure and for this, we take the limit of the stochastic integral equation. From (3.2), the martingale term is asymptotically trivial. Hence, the work is to identify the limit of the term involving the current in terms of the empirical density field. Since the sites are distance $\varepsilon$ apart,

$$\varepsilon^{-1}(\phi(\varepsilon y)-\phi(\varepsilon x)) = (y-x)\cdot\nabla\phi(\varepsilon x) + o(1).$$

Fix a direction $e_0$ and call $J=e_0\cdot\nabla\phi$. Our job is to identify the limit of

(3.8) $$\int_0^T \varepsilon^d \sum_{x\in\mathbb{Z}^d/\varepsilon^{-1}\mathbb{Z}^d} J(\varepsilon x)\varepsilon^{-1}w_{x,x+e_0}\,dt$$

for smooth functions $J$ on the $\mathbb{T}^d$ as

(3.9) $$\int_0^T\int_{\mathbb{T}^d} J(\theta)\sum_e D_{e_0e}(m(\theta,t))\,\partial_e m(\theta,t)\,dt,$$

where $m(\theta,t)$ is the density of the limit of the empirical density $\mu_\varepsilon$. It is relatively easy to see that the latter is well approximated by something of the form

(3.10) $$\int_0^T \varepsilon^d \sum_{x\in\mathbb{Z}^d/\varepsilon^{-1}\mathbb{Z}^d, e'} J(\varepsilon x)D_{e_0e}(\bar\eta_x^{\delta_1\varepsilon^{-1}})(2\delta_2)^{-1}(\bar\eta_{x+\delta_2\varepsilon^{-1}e}^{\delta_1\varepsilon^{-1}} - \bar\eta_{x-\delta_2\varepsilon^{-1}e}^{\delta_1\varepsilon^{-1}})\,dt,$$

where $\bar{\eta}_x^\ell = Av_{y \in \Lambda_x^\ell} \eta_y$ and $\Lambda_x^\ell$ is a cube of side length $\ell$ around $x$, since $m(t, \theta) d\theta$ is the weak limit of $\mu_\varepsilon(t, d\theta)$ and $\bar{\eta}_x^{\delta \varepsilon^{-1}}(t) = \mu_\varepsilon(t, \{\theta : |\theta - x| \leq \delta\})$.

In order to make this replacement, one could try to replace (3.8) by something like

$$(3.11) \qquad \int_0^T \varepsilon^d \sum_{x \in \mathbb{Z}^d/\varepsilon^{-1}\mathbb{Z}^d} J(\varepsilon x) \varepsilon^{-1} \tau_x Av_{y \in \Lambda_K} w_{y, y+e_0} \, dt,$$

where $\Lambda_K$ is a box of side length $K$ about the origin. Note that it is not hard to make such a replacement. Performing a summation by parts on the difference between (3.8) and (3.11), one obtains an error of $\int_0^T \Gamma \, dt$, where

$$(3.12) \qquad \Gamma = \varepsilon^d \sum_{x \in \mathbb{Z}^d/\varepsilon^{-1}\mathbb{Z}^d} [J(\varepsilon x) - Av_{y \in \Lambda_K} J(\varepsilon(x+y))] \varepsilon^{-1} w_{x, x+e_0}.$$

One easily estimates from (3.7)

$$(3.13) \qquad \langle \Gamma f \rangle \leq CK\epsilon \|\nabla J\|_\infty Av_{x \in \mathbb{Z}^d/\varepsilon^{-1}\mathbb{Z}^d} \varepsilon^{-1} \langle (\nabla_{x, x+e_0} \sqrt{f})^2 \rangle^{1/2},$$

so that

$$(3.14) \qquad \langle \Gamma f \rangle - \varepsilon^{d-2} D(\sqrt{f}) \leq C' K^2 \varepsilon^2$$

for some new $C' < \infty$. So as long as $K = o(\varepsilon^{-1})$, such a replacement can be performed. On the other hand, it is not so clear how replacing the current by its average really helps us to get closer to something like (3.10) or where the nontrivial term $\sigma(m)$ would come from.

Instead we will choose some functions $\Theta_K^e$ which depend on variables in box $\Lambda_K$ of side length $K$, large with $\varepsilon^{-1}$, so that the latter is well approximated by

$$(3.15) \qquad \int_0^T \varepsilon^d \sum_{x \in \mathbb{Z}^d/\varepsilon^{-1}\mathbb{Z}^d, e} J(\varepsilon x) \varepsilon^{-1} \tau_x \Theta_K^e \, dt.$$

The functions $\Theta_K^e$ are given explicitly [see (3.28)], but for now we do not need the explicit form to explain the basic argument.

There are also a class of objects for which one can readily check the asymptotics are trivial. Let $g(\eta, \alpha)$ be any local function and define the shift $\tau_x g = g(\tau_x \eta, \tau_x \alpha)$. Then

$$(3.16) \qquad \int_0^T \varepsilon^d \sum_{x \in \mathbb{Z}^d/\varepsilon^{-1}\mathbb{Z}^d} J(\varepsilon x) \varepsilon^{-1} L \tau_x g \, dt$$

is asymptotically trivial. This can be seen by Itô's formula, which says that the above is equal to

$$\varepsilon^{d+1} \sum_x J(\varepsilon x) \tau_x g \Big|_0^T - \int_0^T \varepsilon^{d+1} \sum_{x, e} J(\varepsilon x) a_{x, x+e} \nabla_{x, x+e} g \, dM_{x, x+e}.$$

Terms such as (3.16) are called *fluctuation terms*. Let

$$(3.17) \qquad V_K = Av_{x \in \Lambda'_K} w_{x,x+e_0} - L\tau_x g - \sum_e \Theta^e_K,$$

where $\Lambda'_K$ is a box of side length $K - \ell$, where $\ell$ is chosen so that $V_K$ depends only on variables in $\Lambda'_K$. We need to show that

$$(3.18) \qquad \int_0^T \varepsilon^d \sum_{x \in \mathbb{Z}^d/\varepsilon^{-1}\mathbb{Z}^d} J(\varepsilon x)\varepsilon^{-1} \tau_x V_K \, ds \to 0 \qquad \text{in } P_\varepsilon \text{ probability.}$$

From either the entropy bound in (3.3) together with the entropy inequality, or directly from the Dirichlet form bound in (3.3) to show (3.18) for a function $V_K$, it suffices to prove

$$(3.19) \qquad \inf_g \limsup_{\varepsilon \to 0} \sup_f \left\{ \varepsilon^d \sum_x \varepsilon^{-1} \langle \tau_x V_K f \rangle - \varepsilon^{-2} \sum_e \langle (\nabla_{x,x+e} \sqrt{f})^2 \rangle \right\} = 0.$$

Taking the sum out of the supremum, we obtain an upper bound of the form

$$\inf_g \limsup_{\varepsilon \to 0} \sup \left\{ \varepsilon^{-1} \langle Av_{x \in \Lambda_K} V_K f \rangle - \varepsilon^{-2} Av_{x \in \Lambda_K} \sum_e \langle (\nabla_{x,x+e} \sqrt{f})^2 \rangle \right\}.$$

The expectation is now over a canonical measure (fixed density of particles) on the box $\Lambda^K$ of side length $K$ and the supremum is over all relative density functions, as well as all densities. Letting $L_K$ denote the generator corresponding to the Dirichlet form $\sum_{x \in \Lambda_K} \sum_e \langle (\nabla_{x,x+e} \sqrt{f})^2 \rangle$, we recognize this as the variational formula for

$$\varepsilon^{-2} K^{-(d+2)} \lambda_{K,\varepsilon},$$

where $\lambda_{K,\varepsilon}$ is the principle eigenvalue of

$$K^2 L_K + \varepsilon K^{d+2} V_K.$$

Note that $K^2 L_K$ is used because it has a spectral gap of order one. If we write down the formal Rayleigh–Schrödinger series for $\lambda_{K,\varepsilon}$, we find

$$\varepsilon^{-2} K^{-(d+2)} \lambda_{K,\varepsilon}$$
$$= \varepsilon^{-2} K^{-(d+2)} \{ \varepsilon K^{d+2} \langle V_K \rangle + \varepsilon^2 K^{2(d+2)} \langle V_K, (-K^2 L_K)^{-1} V_K \rangle + \cdots \}.$$

The functions $\Theta^e_K$ are specially chosen so that they have mean zero with respect to any canonical measure on $\Lambda_K$. This is also true of the currents and the fluctuation terms and, hence, $\langle V_K \rangle = 0$. Therefore,

$$(3.20) \qquad \int_0^T \varepsilon^d \sum_{x \in \mathbb{Z}^d/\varepsilon^{-1}\mathbb{Z}^d} J(\varepsilon x)\tau_x [\varepsilon^{-1} V_K - a_K] \, ds$$

is asymptotically trivial where

$$(3.21) \qquad a_K = K^{d+2} \langle V_K, (-K^2 L_K)^{-1} V_K \rangle.$$

In Sections 4 and 5 this is proved for $K \leq c_0 \varepsilon^{-2/(d+2)}$, where $c_0$ is a small constant as long as we have an estimate of the form

$$(3.22) \qquad \langle V_K f \rangle^2 \leq C K^d \sum_{|x-y|=1, x,y \in \Lambda_K} \langle (\nabla_{xy} \sqrt{f})^2 \rangle$$

holding for some $C < \infty$ for all densities $f$ on $\Lambda_K$. For our special choice of $\Theta_K^e$ [see (3.28)], (3.22) will follow from (3.7) and the moving particles lemma proved in Section 5. We want $K$ as large as possible to control fluctuations from the random field and, hence, we will always choose

$$(3.23) \qquad K = c_0 \varepsilon^{-2/(d+2)}.$$

Let

$$(3.24) \qquad \hat{a}_K = a_K - \mathbb{E}[a_K].$$

We claim that, for any fixed $g$, almost surely in the random field,

$$(3.25) \qquad \lim_{\varepsilon \to 0} \varepsilon^d \sum_{x \in \mathbb{Z}^d/\varepsilon^{-1}\mathbb{Z}^d} J(\varepsilon x) \tau_x \hat{a}_K = 0$$

in $P_\varepsilon$ probability. To see this, note that $\hat{a}_K$, which is a function of $m = \bar{\eta}_K$ taking values in $iK^{-d}$, $i = 0, \ldots, K^d$, can be easily extended to a continuous function of $m \in [0,1]$ by linear interpolation. For each $\delta$ with $\delta^{-1}$ a positive integer, divide $\mathbb{Z}^d/\varepsilon^{-1}\mathbb{Z}^d$ into disjoint boxes of side length $\delta \varepsilon^{-1}$ and label them $\beta$. Let $\bar{\eta}_\beta$ be the particle density in box $\beta$. By the two block estimate,

$$(3.26) \qquad \lim_{\delta \to 0} \lim_{\varepsilon \to 0} \varepsilon^d \operatorname*{Av}_{\beta} Av_{x \in \beta} J(\varepsilon x) [\hat{a}_K(\bar{\eta}_{x,K}, \tau_x \alpha) - \hat{a}_K(\bar{\eta}_\beta, \tau_x \alpha)] = 0$$

in $P_\varepsilon$ probability. Now $\hat{a}_K(m)$ is a function of $\alpha_x$, $x \in \Lambda_K$, hence, by the mixing condition (2.1), for $\gamma > 1$,

$$\mathbb{E}\left[\left|\operatorname*{Av}_{\beta} Av_{x \in \beta} J(\varepsilon x) \hat{a}_K(\bar{\eta}_\beta, \tau_x \alpha)\right|^\gamma\right] \leq \operatorname*{Av}_{\beta} \mathbb{E}[|Av_{x \in \beta} J(\varepsilon x) \hat{a}_K(\bar{\eta}_\beta, \tau_x \alpha)|^\gamma]$$

$$= O((\delta \varepsilon K)^{\gamma d/2}).$$

For fixed $\delta$, this is $O(\varepsilon^{\gamma d^2/2(d+2)})$ which is summable in $\varepsilon^{-1} = 1, 2, \ldots$ as long as $\gamma > 2(d+2)/d^2$. By Chebyshev's inequality and the Borel–Cantelli lemma, the term goes to zero for almost every realization of the random field. This proves (3.25).

Hence, we have reduced the problem to proving that

$$(3.27) \qquad \inf_g \lim_{K \to \infty} \sup_m \mathbb{E}[a_K(m, g)] = 0$$

for our specific choice of $\Theta_K^e$, as well as proving that (3.10) is well approximated by (3.15).

The choice of functions $\Theta_K^e$ is not unique. The standard choice is something of the form $D\nabla_e\eta$, where $\nabla_e\eta = \eta_e - \eta_0$, making the passage from (3.15) to (3.10) easy. However, $\nabla_e\eta$ does not have mean zero, so one has to subtract a term $E[\nabla_e\eta|\bar\eta_l]$ for some large $l$, and try to deal with that term in a different way. One can check the size of the subtracted term after appropriate averaging is only small in dimensions three or higher. In [17] it was suggested that the standard approach could work in $d \geq 3$ and it was carried out in [6].

We choose instead in (3.15)

$$\Theta_K^e = \nu(\bar\eta_{\Lambda_K}) W_K^e, \tag{3.28}$$

where, for any integer $\ell$,

$$W_\ell^e = \ell^{-1} \bar w_{\Lambda^\ell, \tau_{\ell e}\Lambda^\ell} \tag{3.29}$$

is the block renormalized or long jump current. Here $\Lambda^\ell = \{1,\ldots,\ell\}^d$ and for any two nonintersecting subsets $A$ and $B$ of $\mathbb{Z}^d$,

$$\bar w_{A,B} = Av_{x\in A, y\in B} w_{xy} \tag{3.30}$$

is the average current over $A$ and $B$, where $w_{xy}$ are given by (3.1). The prefactor $\nu$ is given by

$$\nu(m) = \sigma_{e_0 e}(m)/m(1-m). \tag{3.31}$$

Computations in Section 8 explain why $\nu$ has to have this form. (3.27) is proved in Sections 6 and 7.

**4. Perturbation theory.** We recall some standard results from perturbation theory which we will be using in a specific context.

LEMMA 4.1. *Let $H_0$ be a nonnegative self-adjoint operator on a Hilbert space with $\lambda_0 = 0$ a simple eigenvalue, and spectral gap $\lambda_1 \geq 1$. Let $V$ be a real potential bounded by*

$$\langle u, Vu\rangle \leq \alpha(\langle u, H_0 u\rangle + 1/4\langle u, u\rangle) \tag{4.1}$$

*for some $\alpha \leq 2/5$. Let $H = H_0 + V$. Then for $|\lambda| = 1/2$, the Green function $(\lambda - H)^{-1}$ has the convergent series*

$$(\lambda - H)^{-1} = \sum_{n=0}^\infty [(\lambda - H_0)^{-1} V]^n (\lambda - H_0)^{-1}, \tag{4.2}$$

*with*

$$\|[(\lambda - H_0)^{-1} V]^n (\lambda - H_0)^{-1}\| \leq 5(5\alpha/2)^n. \tag{4.3}$$

PROOF. Note that $(H_0 + 1/4)^{-1/2}$ is well defined. By (4.1),

$$(4.4) \quad \langle u, (H_0 + 1/4)^{-1/2} V (H_0 + 1/4)^{-1/2} u \rangle \leq \alpha \|u\|_2^2.$$

Since $V$ is real, $(H_0 + 1/4)^{-1/2} V (H_0 + 1/4)^{-1/2}$ is self-adjoint, and therefore (4.4) implies the operator bound

$$(4.5) \quad \|(H_0 + 1/4)^{-1/2} V (H_0 + 1/4)^{-1/2}\| \leq \alpha.$$

Rewrite $[(\lambda - H_0)^{-1} V]^n$ as

$$(4.6) \quad \begin{aligned} & (H_0 + \tfrac{1}{4})^{-1/2} [(H_0 + \tfrac{1}{4})(\lambda - H_0)^{-1} (H_0 + \tfrac{1}{4})^{-1/2} V (H_0 + \tfrac{1}{4})^{-1/2}]^n \\ & \times (H_0 + \tfrac{1}{4})^{1/2}. \end{aligned}$$

Since $|\lambda| = 1/2$ and $H_0$ has gap greater than one, we have $\|(H_0 + 1/4)^{-1/2}\| \leq 2$, $\|(H_0 + 1/4)(\lambda - H_0)^{-1}\| \leq 5/2$, $\|(H_0 + 1/4)^{1/2} (\lambda - H_0)^{-1}\| \leq 5/2$. This proves the bound and the convergence. □

The spectral projection for the ground state of $H$ is given by

$$(4.7) \quad P = \frac{1}{2\pi i} \int_{|\lambda| = 1/2} (\lambda - H)^{-1}.$$

Hence, we have a convergent expansion for $P$. Using this expression, one obtains the familiar Rayleigh–Schrödinger series (see [19]) for the ground state energy of $H + V$ with the $n$th term bounded by $5(5\alpha/2)^n$.

COROLLARY 4.2. *Let $W$ be a real potential and*

$$(4.8) \quad H = -K^2 L_K + \epsilon K^{d+2} W.$$

*If $K^2 L_K$ has gap of order one, then the Rayleigh–Schrödinger series for the ground state energy converges provided $W$ satisfies*

$$(4.9) \quad \langle u, W u \rangle \leq K^{-d/2} D_K(u)^{1/2} \|u\|_2, \quad K < c_0 \epsilon^{-2/(d+2)}$$

*for some $c_0$ small enough, or*

$$(4.10) \quad \|W\|_\infty \leq C_1, \quad K \leq (C_1 \epsilon/10)^{-1/(d+2)}.$$

*Furthermore, in both cases one has*

$$(4.11) \quad \epsilon^{-2} K^{-d-2} \operatorname{infspec} H \leq K^d \langle W, -(L_K)^{-1} W \rangle + o(1).$$

PROOF. Assume the gap of $K^2 L_K$ is one. From (4.9) and Schwarz's inequality,

$$(4.12) \quad \epsilon K^{d+2} \langle u, W u \rangle \leq \epsilon K^{(d+2)/2} [K^2 D_K(u) + 1/4 \|u\|_2^2].$$

By Lemma 4.1, the power series for the Green function converges. Furthermore, the $n$th term is bounded by $c(\epsilon K^{(d+2)/2})^n$. Therefore, the $n$th term in the Rayleigh–Schrödinger series is bounded by $c(\epsilon K^{(d+2)/2})^n$. One can compute that the first term is zero from the assumption (4.9) on $W$. The second term gives

$$(4.13) \qquad (\epsilon^2 K^{d+2}) K^d \langle W, -(L_K)^{-1} W \rangle,$$

while the other terms are bounded by

$$(4.14) \qquad (\epsilon K^{(d+2)/2})^3.$$

Now suppose that (4.10) holds instead. We can choose $\alpha$ in the previous lemma to be $\epsilon K^{d+2} C_1$. The first two terms in the Rayleigh–Schrödinger series can be computed as before. The other terms are again bounded by $o(1)$ after multiplying $\epsilon^{-2} K^{-d-2}$. This proves the corollary. $\square$

**5. Moving particle lemma.** We need a preliminary result.

LEMMA 5.1. *Suppose $\mu$ is a homogeneous Bernoulli measure ($\alpha_x \equiv 0$) on $\mathbb{Z}^1$ conditioned to have $N$ particles. Let $k$ be a positive integer and $\rho_x$ a sequence of positive numbers with $\sum_{x=1}^{k-1} \rho_x = 1$. Then*

$$(5.1) \qquad E_\mu[(\nabla_{1,k} f)^2] \le \sum_{x=1}^{k-1} \rho_x^{-1} E_\mu[(\nabla_{x,x+1} f)^2].$$

It is important that there is no multiplicative constant on the right-hand side.

PROOF OF LEMMA 5.1. We prove it by induction. For $k=3$, it is elementary to check the inequality directly. Suppose it is true for $k-1$ and let $\rho_1 + \cdots + \rho_k = 1$. Let $q = \rho_2 + \cdots + \rho_k$. From the $k=3$ case we have

$$(5.2) \qquad E[(\nabla_{1,k} f)^2] \le \rho_1^{-1} E[(\nabla_{1,2} f)^2] + q^{-1} E[(\nabla_{2,k} f)^2].$$

From the inductive hypothesis, since $q^{-1}\rho_2 + \cdots + q^{-1}\rho_k = 1$,

$$(5.3) \qquad E[(\nabla_{2,k} f)^2] \le q \sum_{x=2}^{k-1} \rho_x^{-1} E_\mu[(\nabla_{x,x+1} f)^2].$$

This completes the proof. $\square$

LEMMA 5.2 (Moving particle lemma). *Suppose $\mu$ is a Bernoulli measure on $\mathbb{Z}^1$ with external field $\alpha$ taking values in $[-B, +B]$. Then*

$$(5.4) \qquad \int [f(T_{1L}\eta) - f(\eta)]^2 \, d\mu(\eta) \le e^{13B} L \sum_{1 \le x \le L-1} \int (\nabla_{x,x+1} f)^2 \, d\mu(\eta).$$

PROOF. Suppose we change the measure $\mu$ to a new measure $\tilde{\mu}$ by changing each $\alpha_i$ to the nearest value of the form $Bj/L$, $j$ an integer. Since there is always such a point with $|\alpha_i - Bj/L| \leq B(2L)^{-1}$, the Radon–Nikodym derivative $d\mu/d\tilde{\mu}$ is bounded above and below uniformly by $e^{B/2}$ and $e^{-B/2}$, respectively. Therefore, at the cost of a factor of $e^B$, we may assume that $\alpha$ takes values in $\{Bj/L : j = -L, \ldots, L\}$. By the same reasoning, at the price of a factor $e^{4B}$, we may assume that $\alpha_0 = \alpha_L = B$.

Let $A$ be the set
$$A = \{x_i : \alpha_{x_i} = K, i = 1, \ldots, k\}.$$

By definition,
$$T_{1,L}\eta = T_{x_1,x_2} \cdots T_{x_{k-2},x_{k-1}} T_{x_k,x_{k-1}} \cdots T_{x_3,x_2} T_{x_2,x_1} \eta.$$

By Lemma 5.1 with $\rho_s^{-1} = \frac{L}{x_{s+1}-x_s}$,

$$E_\mu[(f(T_{1,L}\eta) - f(\eta))^2] \leq e^{4K} \sum_{s=1}^{k-1} \frac{L}{x_{s+1} - x_s} E_\mu[(f(T_{x_s,x_{s+1}}\eta) - f(\eta))^2].$$

We have to bound
$$\frac{L}{x_{s+1} - x_s} E_\mu[(f(T_{x_s,x_{s+1}}\eta) - f(\eta))^2].$$

We are now in the same situation as before except no $\alpha_x$ can take value $K$ when $x_s < x < x_{s+1}$. Let us change $\alpha_{x_s}$ and $\alpha_{x_{s+1}}$ to the value $B(L-1)/L$. The price we pay is a factor $\exp\{2BL^{-1}\}$. Continuing this procedure, we have a proof of the lemma. □

From the moving particles lemma, we obtain as a consequence two basic bounds, the spectral gap and the two block estimate. The arguments leading from Lemma 5.2 to the two block estimate are completely standard (see [11]), so we simply state the result we need. We denote by $\bar{\eta}_x^n$ the empirical density of particles on a box of side length $n$ around $x$. Also $\langle \cdot \rangle_m$ denotes the expectation with respect to the fixed total density of particles $\bar{\eta}_{\varepsilon^{-1}} = m$.

LEMMA 5.3 (Two block estimate). *Let $0 \leq K < \infty$ and $\alpha$ any field with $-K \leq \alpha_x \leq K$. Let $F$ be continuous function on $[0,1] \times [0,1]$ and $J$ a smooth test function on $\mathbb{T}^d$. Let $K \to \infty$ as $\varepsilon \to 0$ with $K \leq \delta \varepsilon^{-1}$. Let*

(5.5) $$V_{\varepsilon,\delta} = \varepsilon^d \sum_x J(\varepsilon x)(F(\bar{\eta}_{x+y}^K, \bar{\eta}_{x+z}^K) - F(\bar{\eta}_x^{\delta\varepsilon^{-1}}, \bar{\eta}_x^{\delta\varepsilon^{-1}})).$$

*Then for every realization of the field $\alpha$ with $-K \leq \alpha_x \leq K$,*

$$\limsup_{\delta \to 0} \limsup_{\varepsilon \to \infty} \sup_{\substack{|y|,|z| \leq K \\ 0 \leq m \leq 1}} \sup_f \left\{ \langle V_{\varepsilon,\delta} f \rangle_m \right.$$

(5.6)

$$\left. - \varepsilon^{d-2} \sum_{x \in \mathbb{Z}^d/\varepsilon^{-1}\mathbb{Z}^d, e} \langle (\nabla_{x,x+e}\sqrt{f})^2 \rangle_m \right\} \leq 0.$$

*The expectation $\langle \cdot \rangle_m$ is with respect to the canonical measure (2.8) with fixed density $m$ on $\mathbb{Z}^d/\varepsilon^{-1}\mathbb{Z}^d$ and the supremum is over all relative densities $f$.*

The following result from [4] gives the spectral gap of the Bernoulli–Laplace version of our model to the correct order.

LEMMA 5.4. *Let $0 \leq K < \infty$. There exists a constant $C = C(K) < \infty$ such that, for any field $\alpha$ with $-K \leq \alpha_x \leq K$, any $\Lambda$, and any $0 \leq N \leq |\Lambda|$, for any $f : \{0,1\}^\Lambda \to \mathbb{R}$,*

(5.7) $$\operatorname{Var}_{\Lambda,N}(f) \leq C|\Lambda|^{-1} E_{\Lambda,N}\left[ \sum_{x,y \in \Lambda} (\nabla_{xy} f)^2 \right].$$

*Here $\operatorname{Var}_{\Lambda,N,K}$ and $E_{\Lambda,N,K}$ denote the variance and expectation with respect to the measure in our random field on $\Lambda$ conditioned to have $N$ particles.*

Together with the moving particles lemma, one obtains the spectral gap of the nearest neighbor dynamics to the correct order.

THEOREM 2 (Spectral gap). *For each $K > 0$, there exists a $C < \infty$ such that, for all $\alpha$ with $-K \leq \alpha_x \leq K$, all cubes $\Lambda_L$ of side length $L$, all $0 \leq N \leq L^d$, and all $f : \{0,1\}^{\Lambda_L} \to \mathbb{R}$,*

(5.8) $$\operatorname{Var}_{\Lambda_L,N}(f) \leq CL^2 E_{\Lambda_L,N}\left[ \sum_{\substack{x,y \in \Lambda_L \\ |x-y|=1}} (\nabla_{xy} f)^2 \right].$$

PROOF. By Lemma 5.4, we have

$$\operatorname{Var}_{\Lambda_L,N}(f) \leq \frac{C}{|\Lambda_L|} E_{\Lambda_L,N}\left[ \sum_{x,y \in \Lambda_L} (\nabla_{xy} f)^2 \right].$$

For each $x, y \in \Lambda_L$, choose a canonical path $x = x_1, x_2, \ldots, x_n = y$ with $x_i$ and $x_{i+1}$ by moving first in the first coordinate direction, then in the second coordinate direction, and so on. By the moving particles lemma, we have

$$E_{\Lambda_L,N}[(f(T_{xy}\eta) - f(\eta))^2] \leq e^{13K} n \sum_{1 \leq i \leq n-1} E_{\Lambda_L,N}[(\nabla_{x_i,x_{i+1}} f)^2].$$

Summing over $x$ and $y$, noting that $n \leq dL$ and that each nearest neighbor pair is used for the path between $d(L/2)^{d+1}$ pairs $x$ and $y$, we obtain the result. □

**6. The diffusion coefficient.** We now discuss the computation of the asymptotic variance (3.27). Suppose $f(\eta, \alpha)$ is a function depending on the particle and field configuration in some finite box $\Lambda_\ell$ and which has mean zero with respect to every invariant measure for our process on that box. We denote the class of such functions $\mathcal{G}_\ell$. Denote by $\mathcal{G}$ the increasing union of $\mathcal{G}_\ell$. For any $v \in \mathcal{G}_\ell$, we form the sum

$$V_K = \sum_{x \in \Lambda'_K} \tau_x v, \tag{6.1}$$

where $\Lambda'_K$ is a box of side length $K - \ell$ and we define

$$\mathbf{V}_K = \langle V_K, (-L_K)^{-1} V_K \rangle = \sup_f \left\{ 2 \langle V_K f \rangle - \sum_{\substack{|x-y|=1 \\ x,y \in \Lambda_K}} \langle (\nabla_{xy} f)^2 \rangle \right\}. \tag{6.2}$$

The expectation is with respect to an extremal invariant measure, that is, the canonical invariant measure (2.8) with fixed number of particles, and the supremum is over all densities with respect to such a measure. In particular, $\mathbf{V}_K$ depends on $\bar{\eta}_{\Lambda_K} = Av_{x \in \Lambda_K} \eta_x$, as well as $\alpha_x$, $x \in \Lambda_K$,

$$\mathbf{V}_K = \mathbf{V}_K(\bar{\eta}_{\Lambda_K}, \alpha). \tag{6.3}$$

Since $v \in \mathcal{G}_\ell$, if $K \gg \ell$, we can write $v = L_K h$ for some $h$ supported on $\Lambda_K$. Then $V_K = \sum_{x \in \Lambda'_K} \tau_x L_K h = L_K \sum_{x \in \Lambda'_K} \tau_x h$ and, hence, we have $\mathbf{V}_K \leq CK^d$. For each fixed $m \in (0,1)$, define

$$[v, v] = \limsup_{m_K \to m} K^{-d} \mathbb{E}[\mathbf{V}_K(m_K, \alpha)]. \tag{6.4}$$

Although $[v, v]$ is well defined for each $v \in \mathcal{G}$, it is not so clear how to compute it. For this purpose, we introduce an auxiliary Hilbert space $\mathcal{H}$ and compute $[v, v]$ for objects in $\mathcal{G}$ by mapping them to $\mathcal{H}$. Let $\pi_m$ denote the product measure (2.6) with density $m$. A *form* $\xi_e$ is a function of $\eta$ and $\alpha$ depending on the basic basis elements $e$ of $\mathbb{Z}^d$. We make a Hilbert space $\mathcal{H}$ out of them through the inner product

$$\langle\!\langle \xi, \xi \rangle\!\rangle = \sum_e \mathbb{E}[\langle \xi_e^2 \rangle_m]. \tag{6.5}$$

For $g$ a local function, consider the *exact* form $\xi_b = \sum_{y \in \mathbb{Z}^d} \nabla_b \tau_y g$ indexed over nearest neighbor bonds $b = x, x+e$ in $\mathbb{Z}^d$. Since $g$ is local, there are only finitely many nonzero terms so the sum is well defined. Note that $\xi_b$ is covariant in the sense that $\tau_x \xi_b = \xi_{\tau_x b}$. So they can all be reconstructed out of the basic forms $\xi_e$, $e$ running over basis elements of $\mathbb{Z}^d$. This defines a Hilbert subspace $\mathcal{E}$ of exact covariant forms. An exact form possesses the algebraic property of *closedness: If $b_1, b_2, \ldots, b_n$ are an ordered set of bonds making a loop in $\mathbb{Z}^d$, then $\sum_{i=1}^n \xi_{b_i}(\eta^{b_1 \cdots b_{i-1}}) = 0$.* An example of a form which is closed but not exact is given by $(\nabla \eta)_e = \nabla_e \eta = \eta_e - \eta_0$. Let $\mathcal{C}$ be the closure in $\mathcal{H}$ of the closed forms. One can check by standard arguments (see [11]) that $\mathcal{E}$ has codimension $d$ in $\mathcal{C}$ and that $\mathcal{C}$ is the closure of the linear span of $\mathcal{E}$ and $\nabla \eta$.

We now compute $[v, v]$ for some special cases. When $v = Lg$ for some local $g$, we use the map

$$(6.6) \qquad g \mapsto \xi_e = \nabla_{0,e} \sum_x \tau_x g$$

and find

$$(6.7) \qquad [Lg, Lg] = \langle\!\langle \xi, \xi \rangle\!\rangle.$$

When $v = w_{0,e}$, we have the image $w \mapsto \nabla \eta$ and

$$(6.8) \qquad [w, w] = \langle\!\langle \nabla \eta, \nabla \eta \rangle\!\rangle.$$

Furthermore, for $v = Lg$, the map above gives

$$(6.9) \qquad [w_{0,e}, v] = -\tfrac{1}{2} \mathbb{E}[\langle \nabla_e \eta, \xi_e \rangle_m].$$

The variational formula (2.12) for $\sigma$ translates to

$$(6.10) \qquad \beta \cdot \sigma \beta = \inf_g \left[ \left( \sum_e \beta_e w_{0,e} - Lg \right)^2 \right].$$

Fix a direction $e$. For each fixed $\ell$, $W_\ell^e$ given by (3.29) is embedded in $\mathcal{E}$, as $\varpi_\ell^e$. From (3.7) and the moving particles lemma, one has a uniform bound on the norm $[W_\ell^e, W_\ell^e]$. Furthermore, one can check by a straightforward computation that, for any $g \in \mathcal{G}$,

$$(6.11) \qquad \lim_{\ell \to \infty} [W_\ell^e, Lg] = 0$$

and less trivially that

$$(6.12) \qquad \lim_{\ell \to \infty} [\beta \cdot W_\ell, \beta \cdot W_\ell] = m^2 (1-m)^2 \beta \cdot \sigma^{-1} \beta.$$

In other words, the limit $W$ of the $W_\ell$ represents an (unnormalized) orthogonal complement to $Lg$ with respect to $[\cdot, \cdot]$. From the above computations, we have

$$(6.13) \qquad \lim_{\ell \to \infty} \inf_g [(w - Lg - \nu W_\ell)^2] = 0,$$

where $\nu(m) = \sigma(m)/m(1-m)$, which is the first step in the renormalization.

To prove (6.12), we take the limit in variational formula (6.2) (see [11]) for any $v \in \mathcal{G}_\ell$ to find

$$(6.14) \qquad [v,v] = \sup_\xi \left\{ 2 \sum_{x, x+e \in \Lambda_\ell} \mathbb{E}[\langle \nabla_e \tau_{-x} V, \tau_{-x} \xi_e \rangle] - \langle\langle \xi, \xi \rangle\rangle \right\},$$

where $-L_{\lambda_\ell} V = v$ and $\xi$ is a closed form. Since $\xi$ is closed, it can be approximated by $\beta \nabla \eta + \nabla \sum_y \tau_y g$. The first term can be computed explicitly and the variational formula becomes

$$(6.15) \qquad [v,v] = \sup_{g \in \mathcal{G}, \beta} \left\{ 2 \mathbb{E}\left[ \left\langle v, \sum_y \tau_y g + \sum \beta \cdot y \eta_y \right\rangle \right] \right.$$
$$\left. - \left\langle\!\left\langle \left( \beta \nabla \eta - \nabla \sum_x \tau_x g \right)^2 \right\rangle\!\right\rangle \right\}.$$

Applying this to $W_\ell$ gives (6.12). It is worth noting that

$$\mathbb{E}\left[ \left\langle v, \sum_y \tau_y g + \sum \beta \cdot y \eta_y \right\rangle \right] = [v, \beta \cdot w - Lg]$$

and

$$\left\langle\!\left\langle \left( \beta \nabla \eta - \nabla \sum_x \tau_x g \right)^2 \right\rangle\!\right\rangle = [(\beta \cdot w - Lg)^2],$$

so this also shows that $\mathcal{G}$ is a Hilbert space and $\mathcal{G} = \bar{w} \oplus L\mathcal{G}$.

At this point it is worth making a quick remark about associated dynamics. If one chooses instead a dynamics as in (2.16), then one will have a corresponding generator $\tilde{L}$ and a corresponding current $\tilde{w}$ different from the $w$ above, but satisfying all the needed estimates. One checks in exactly the same way as above that $\lim_{\ell \to \infty} \inf_g [(\tilde{w} - \tilde{L}g - \tilde{\nu} W_\ell)^2] = 0$, where $\tilde{\nu}(m) = \tilde{\sigma}(m)/m(1-m)$ and $\tilde{\sigma}$ is given by (2.17). The key point is that the choice $W_\ell$ is independent of the particular choice of model. This is important in later sections where the special form of $W_\ell$ is used repeatedly.

**7. Structure of the gradient space.** The state space on which we work is $\Omega = [-B, B]^{\mathbb{Z}^d} \times \{0,1\}^{\mathbb{Z}^d}$ points of which we call $(\alpha, \eta)$. Fix $0 < m < 1$. On $\{0,1\}^{\mathbb{Z}^d}$, we have the product measure

$$(7.1) \qquad \mu = Z^{-1} \exp\left\{ \sum_{x \in \varepsilon \mathbb{Z}^d / \mathbb{Z}^d} (\alpha_x + \lambda(m)) \eta_x \right\},$$

with $\lambda(m)$ chosen as in (2.7). Denote corresponding expectation by $E_\alpha[\cdot]$. On $\Omega$ we have the corresponding annealed measure

$$\mathbf{E}[F(\alpha,\eta)] = \mathbb{E}[E_\alpha[F]]. \tag{7.2}$$

A collection of functions $\{\omega_b\}$ where $b$ runs over unoriented bonds of $\mathbb{Z}^d$ is called a *translation covariant form* if, for all bonds $b$ and sites $x$ in $\mathbb{Z}^d$,

$$\tau_x \omega_b = \omega_{\tau_x b}. \tag{7.3}$$

Here the shift operator $\tau_x$ acts on functions $F$ by $(\tau_x F)(\alpha, \eta) = F(\tau_x \alpha, \tau_x \eta)$, where $(\tau_x \eta)_y = \eta_{x+y}$ and $(\tau_x \alpha)_y = \alpha_{x+y}$. A translation covariant form is represented by $\{w_e\}$, where $e$ runs over the basic bonds $e_1 = (1,0,\ldots,0), \ldots, e_d = (0,\ldots,0,1)$ since every $b = \tau_x e$ for some $x$ and $e$ and $\omega_b = \tau_x \omega_e$. The collection of translation covariant forms is then a Hilbert space $\mathcal{H}$ with the norm $\|\{\omega_b\}\| = \sum_e \mathbf{E}[\omega_e^2]$.

An ordered finite set of bonds $b_1, \ldots, b_n$ is called a *closed loop* if $T_{b_n} \cdots T_{b_1} \eta = \eta$ for all $\eta$. A form $\{\omega_b\}$ is called *closed* if for any closed loop of bonds $b_1, \ldots, b_n$,

$$\sum_{i=1}^n \omega_{b_i}(\alpha, T_{b_{i-1}} \cdots T_{b_i} \eta) = 0. \tag{7.4}$$

A translation covariant form is called *exact* if $\omega_e = \sum_x \nabla_e \tau_x g$ for some local function $g(\alpha, \eta)$. Note that only finitely many terms in the sum are finite, so this makes sense. Note also that an exact form is automatically closed.

Denote by $\mathcal{C}$ the closure in $\mathcal{H}$ of the space of square integrable translation covariant closed forms and by $\mathcal{E}$ the closure of those that are exact. We have $\mathcal{E} \subset \mathcal{C}$ and, in fact, the containment is strict as $\nabla_e \eta$, the form $(\eta_{(1,0,\ldots,0)} - \eta_{(0,\ldots,0)}, \ldots, \eta_{(0,\ldots,0,1)} - \eta_{(0,\ldots,0)})$ is closed but not in the closure of the span of exact forms. The main result follows:

THEOREM 3. $\mathcal{C} = \nabla \eta \oplus \mathcal{E}$.

Note that the theorem has been proved in different situations in [6, 17, 22]. The proof here is very similar. It requires that the disorder field $\alpha$ be a stationary ergodic process, but does not use the mixing conditions (2.1).

PROOF OF THEOREM 3. Let $\omega_b$ be a closed form. Let $\Lambda_K$ be a box in $\mathbb{Z}^d$ of side length $K$ centered at the origin. For any set $\Lambda$, let $\eta_\Lambda, \alpha_\Lambda$ denote the collection of variables $\eta_x, \alpha_x, x \in \Lambda$. For bonds $b$ inside $\Lambda_K$, let

$$\omega_b^K = \mathbf{E}[\omega_b | \alpha_{\Lambda_K}, \eta_{\Lambda_K}]. \tag{7.5}$$

$\omega_b^K$ satisfies the closedness condition (7.4) for closed loops $b_1, \ldots, b_n$ contained inside $\Lambda_K$. The set $\{0,1\}^{\Lambda_K}$ is divided into ergodic classes $\sum_{x \in \Lambda_K} \eta_x =$

$N$ and on each we can sum along bonds to produce an unambiguous function $S_K$ with

$$\nabla_b S_K = \omega_b^K \tag{7.6}$$

for any bond $b$ inside $\Lambda_K$. There is a free variable depending on $\sum_{x \in \Lambda_K} \eta_x$ and, therefore, we can choose $S_K$ so that $E_{\mu_m}[S_K | \sum_{x \in \Lambda_K} \eta_x] = 0$ as well. Let

$$G_K = \mathbf{E}[S_{3K} | \alpha_{\Lambda_K}, \eta_{\Lambda_K}] \tag{7.7}$$

and, for any bond $b$ in $\mathbb{Z}^d$,

$$\xi_b^K = K^{-d} \sum_{x \in \mathbb{Z}^d} \nabla_b \tau_x G_K. \tag{7.8}$$

Note that $\xi_b^K$ makes sense because only finitely many terms are nonzero, and that it is by definition an exact form. If $e$ is a basic bond, we can write $\xi_e^K = \xi_e^{1,K} + \xi_e^{2,K}$, where $\xi_e^{1,K}$ corresponds to terms in the sum with $x$ and $x+e$ in $\Lambda_K$ and $\xi_e^{2,K}$ corresponds to terms in the sum with one of $x$ and $x+e$ in $\Lambda_K$ and one not in $\Lambda_K$:

$$\begin{aligned}
\xi_e^{1,K} &= K^{-d} \sum_{x, x+e \in \Lambda_K} \tau_{-x} \nabla_{x,x+e} \mathbf{E}[S_{3K} | \alpha_y, \eta_y, y \in \Lambda_K] \\
&= K^{-d} \sum_{x, x+e \in \Lambda_K} \mathbf{E}[\tau_{-x} \nabla_{x,x+e} S_{3K} | \alpha_{\tau_{-x}\Lambda_K}, \eta_{\tau_{-x}\Lambda_K}] \\
&= K^{-d} \sum_{x, x+e \in \Lambda_K} \mathbf{E}[\omega_e | \alpha_{\tau_{-x}\Lambda_K}, \eta_{\tau_{-x}\Lambda_K}].
\end{aligned} \tag{7.9}$$

By the martingale convergence theorem, we have $\xi_e^{1,K} \to \omega_e$ as $K \to \infty$.

Note that

$$\xi_e^{2,K} = c_e^K \nabla_e \eta, \tag{7.10}$$

where $c_e^K$ does not depend on $\eta_0$ or $\eta_e$. Suppose we take a bond $b$ which does not have any vertex in common with $e$:

$$\nabla_b \xi_e^{2,K} = K^{-d} \sum_{|\{x,x+e\} \cap \Lambda_K|=1} \nabla_e \tau_x \nabla_{\tau_x b} G_K \tag{7.11}$$

so

$$\mathbf{E}[|\nabla_b \xi_e^{2,K}|^2] \leq C K^{-(d+1)} \sum_{|\{x,x+e\} \cap \Lambda_K|=1} \mathbf{E}[|\nabla_{x+b} G_K|^2] \leq C' K^{-2}. \tag{7.12}$$

In other words, $\mathbf{E}[(\nabla_b c_e^K)^2] \to 0$ as $K \to \infty$. Since $c_e^K$ does not depend on $\eta_0$ or $\eta_e$ and this is true for any $b$ not having a vertex in common with $e$, any limit $c_e$ of the $c_e^K$ cannot depend on the configuration $\eta$.

We have thus shown that if $\omega_b$ is a closed form, then there exist exact functions $\xi_e^K$ and a function $\xi_e^{2,K}$ such that $\xi_e^K - \xi_e^{2,K} \to \omega_e$ and any weak limit $\xi_e^2$ of the $\xi_e^{2,K}$ must be of the form

$$c_e(\alpha) \nabla_e \eta. \tag{7.13}$$

Assume for a moment that such a weak limit exists. We now want to show that a form of the type (7.13) is in $\nabla_e \eta \oplus \mathcal{E}$ as well.

First of all, we can always subtract $c\nabla_e \eta$ from (7.13), so taking $c = \mathbb{E}[c_e(\alpha)]$, we can assume without loss of generality that $\mathbb{E}[c_e] = 0$.

Now we solve explicitly for a function $\tilde{G}_K$ such that $\nabla_b \tilde{G}_K = c_e(\tau_x \alpha) \nabla_b \eta$ for bonds $b$ inside $\Lambda_K$ of the form $(x, x+e)$. Trying $\tilde{G}_K = \sum_{x \in \Lambda_K} a_x(\alpha) \eta_x$, one obtains the set of equations

$$a_x - a_y = c_e(\tau_x \alpha). \tag{7.14}$$

These are readily solved by taking $a_0 = 0$ and for any other $x \in \Lambda_K$,

$$a_x = -\sum_{i=1}^{n} c_{e_i}(\tau_{x_i} \alpha), \tag{7.15}$$

where $(x_i, x_i + e_i)$, $i = 1, \ldots, n$, is a sequence of bonds from $0$ to $x$. It does not depend on the path taken because (7.13) is closed.

Let $\tilde{\xi}_e^K$ be as in (7.8) with this new $\tilde{G}_K$, and again break it up as $\tilde{\xi}_e^K = \tilde{\xi}_e^{1,K} + \tilde{\xi}_e^{2,K}$. The first term $\tilde{\xi}_e^{1,K}$ converges to (7.13) in the same way as before. The difference is that we now compute $\tilde{\xi}_e^{2,K}$ directly. It is given by

$$\tilde{\xi}_e^{2,K} = r_e^K \nabla_e \eta, \qquad r_e^K = K^{-d} \sum_{x \in \Lambda, x+e \notin \Lambda_K} \tau_{-x} a_x \tag{7.16}$$

plus an analogous term with $x + e \in \Lambda_K$ and $x \notin \Lambda_K$. Connect boundary points $x$ of $\Lambda_K$ to the origin in some deterministic way: Say, $(x_i(x), x_i(x) + e_i(x))$, $i = 1, \ldots, n(x)$, with $n(x) \leq Kd$. Then

$$r_e^K = -K^{-d} \sum_{x \in \Lambda, x+e \notin \Lambda_K} \sum_{i=1}^{n(x)} c_{e_i(x)}(\tau_{x_i(x)-x} \alpha). \tag{7.17}$$

As $K \to \infty$, this converges to $\mathbb{E}[c_e] = 0$ by the ergodic theorem. This proves that forms which can be written as in (7.13) are in $\nabla_e \eta \oplus \mathcal{E}$.

Finally, we have to prove the key analytic point which is the boundedness in norm of the $\xi_e^{2,K}$ defined after (7.8) and, hence, the existence of a weak limit point. This is in fact the hard point in the nongradient method and the crucial point where the spectral gap is used. We will use the following lemma.

LEMMA 7.1. *Let $\Lambda_1$ and $\Lambda_2$ be disjoint finite subsets of $\mathbb{Z}^d$ and $x$ a point not in either. There exists a constant $C < \infty$ such that, for all functions $f$ of $\eta_{\Lambda_1 \cup \Lambda_2 \cup \{x\}}$, all $\alpha$ and all $y \in \Lambda_2$,*

$$
\begin{aligned}
(7.18) \quad & E_\alpha[(\nabla_{xy} E_\alpha[f|\eta_{\Lambda_1 \cup \{x\}}])^2 | \eta_{\Lambda_1}] \\
& \leq C\{|\Lambda_2|^{-1} E_\alpha[(f - E_\alpha[f|\eta_{\Lambda_1}])^2] + Av_{z \in \Lambda_2} E_\alpha[(\nabla_{xz} f)^2 | \eta_{\Lambda_1}]\}.
\end{aligned}
$$

PROOF. The proof is standard. We write $E[\cdot]$ for $E_\alpha[\cdot|\eta_{\Lambda_1}]$ and $P$ for the corresponding probabilities. We can assume without loss of generality that $E[f] = 0$. Note that

$$(7.19) \quad (\nabla_{xy} E[f|\eta_{\{x\}}])^2 = (E[f|\eta_x = 1] - E[f|\eta_x = 0])^2 \mathbb{1}(\eta_x \neq \eta_y)$$

and

$$
\begin{aligned}
& E[f|\eta_x = 1] - E[f|\eta_x = 0] \\
&= E[f|\eta_x = 1] - E[f|\eta_x = 1, \eta_z = 0] + E[f|\eta_x = 1, \eta_z = 0] \\
(7.20) \quad & \quad - E[f|\eta_x = 0, \eta_z = 1] + E[f|\eta_x = 0, \eta_z = 1] - E[f|\eta_x = 0] \\
&= E\left[f\left(1 - \frac{1-\eta_z}{p_z}\right)\Big|\eta_x = 1\right] + E\left[(\nabla_{xz} f)\frac{1-\eta_z}{p_z}\Big|\eta_x = 1\right] \\
& \quad - E\left[f\left(1 - \frac{\eta_z}{1-p_z}\right)\Big|\eta_x = 0\right],
\end{aligned}
$$

where $p_z = E[\eta_z]$. Since $0 < m < 1$ and $|\alpha_x| \leq B < \infty$, we have a bound $\delta \leq p_z \leq 1 - \delta$ for some $\delta > 0$. By Schwarz's inequality, we therefore have

$$
\begin{aligned}
(7.21) \quad & (E[f|\eta_x = 1] - E[f|\eta_x = 0])^2 \\
& \leq C\{E[f^2|\eta_x = 1] + E[(\nabla_{xz} f)^2|\eta_x = 1] + E[f^2|\eta_x = 0]\}.
\end{aligned}
$$

But again $E[\cdot|\eta_x = 1] \leq C E[\cdot]$ and $E[\cdot|\eta_x = 0] \leq C E[\cdot]$ and the lemma follows. □

Applying the lemma with $f = E_\alpha[S_{3K}|\eta_{\Lambda_K}]$, $x \in \Lambda_K$, $y = x + e \notin \Lambda_K$, $\Lambda_1 = \Lambda_K \setminus \{x\}$, $\Lambda_2 = \{x + e, x + 2e, \ldots, x + Ke\}$, we obtain

$$(7.22) \quad E_\alpha[(\nabla_{x,x+e} f)^2] \leq C\{K^{-d} E_\alpha[S_{3K}; S_{3K}] + Av_{z \in \Lambda_2} E_\alpha[(\nabla_{xz} S_{3K})^2]\}.$$

By the moving particles lemma, Lemma 5.2,

$$
\begin{aligned}
(7.23) \quad E_\alpha[(\nabla_{xz} S_{3K})^2] & \leq CK \sum_{i=1}^{z} E_\alpha[(\nabla_{x+(i-1)e, x+ie} S_{3K})^2] \\
& = CK \sum_{i=1}^{n(z)} E_\alpha[(\omega_{x+(i-1)e, x+ie}^{3K})^2].
\end{aligned}
$$

By the spectral gap, Theorem 2, there is a $C < \infty$ such that, for any $\alpha$,

$$(7.24) \qquad E_\alpha[S_K^2] \leq CK^2 \sum_{b \in \Lambda_K} E_\alpha[|\omega_b^K|^2].$$

Taking expectation of (7.22), (7.23) and (7.24) over $\alpha$, and using Jensen's inequality and $\mathbf{E}[\omega_b] = C < \infty$,

$$(7.25) \qquad \mathbf{E}[(\nabla_{x,x+e} G_K)^2] \leq CK^2$$

and, hence, another application of Jensen's inequality gives the required bound on the norm of $\xi_e^{2,K}$:

$$(7.26) \quad E[|\xi_e^{2,K}|^2] = E\left[\left(K^{-d} \sum_{\text{one of } x,x+e \in \Lambda_K} \tau_{-x} \nabla_{x,x+e} G_K\right)^2\right] \leq C. \qquad \square$$

**8. The long jump current.** So far everything we have described is standard (see [11]). Applying the perturbation theory and the discussion from Section 3, we have obtained that the difference between (3.8) and

$$(8.1) \qquad \int_0^T \varepsilon^d \sum_{x \in \mathbb{Z}^d/\varepsilon^{-1}\mathbb{Z}^d, e} J(\varepsilon x) \tau_x \varepsilon^{-1} \nu_{e_0,e}(\bar{\eta}_{\Lambda_K}) Av_{y \in \Lambda'_K} \tau_y W_\ell^e \, dt$$

goes to 0 in $P_\varepsilon$-probability as $\varepsilon \to 0$ followed by $\ell \to \infty$, where $K = c_0 \varepsilon^{-2/(d+2)}$. However, we will need to take the $\ell$ in $\bar{w}_\ell^e$ on a much larger scale in order for the averaging to beat the fluctuations from the random field. We will again use central limit variance computation as in the previous section, but now we will not be allowed to have the full average over shifts which is crucial in (6.2), so the problem has to be handled in a new way.

To compute the central limit theorem variance of the long jump current, we will use its precise form and a renormalization procedure.

Using its precise form, we can rewrite the current (3.1) as

$$(8.2) \qquad w_{xy} = \eta_y - \eta_x + \xi_x \zeta_y - \zeta_x \xi_y,$$

where

$$(8.3) \qquad \xi_x = (1 - \eta_x) e^{\alpha_x}, \qquad \zeta_x = \eta_x e^{-\alpha_x}.$$

Now we introduce some notation. Let $\Omega$ be a subset of $\mathbb{Z}^d$ containing two nonintersecting subsets $\Lambda_1$ and $\Lambda_2$. The average current on $\Lambda_1, \Lambda_2$ is given by

$$\bar{w}_{\Lambda_1,\Lambda_2} = Av_{x \in \Lambda_1, y \in \Lambda_2} w_{xy}.$$

Let

$$\bar{\eta}_\Lambda = Av_{x \in \Lambda} \eta_x$$

be the empirical density on a set $\Lambda \subset \mathbb{Z}^d$ and $\hat{\lambda}_\Lambda = \hat{\lambda}_\Lambda(m)$ be the empirical chemical potential defined implicitly through

$$(8.4) \qquad m = Av_{x \in \Lambda} \langle \eta_x \rangle_{\hat{\lambda}_\Lambda} = Av_{x \in \Lambda} \frac{e^{\hat{\lambda}_\Lambda + \alpha_x}}{1 + e^{\hat{\lambda}_\Lambda + \alpha_x}}.$$

We will often write $\hat{\lambda}_\Lambda$ for $\hat{\lambda}_\Lambda(\bar{\eta}_\Lambda)$. Define also

$$(8.5) \qquad \bar{\xi}_\Lambda = Av_{x \in \Lambda} \xi_x, \qquad \bar{\zeta}_\Lambda = Av_{x \in \Lambda} \zeta_x.$$

We have

$$(8.6) \qquad \bar{w}_{\Lambda_1, \Lambda_2} = \bar{\eta}_{\Lambda_2} - \bar{\eta}_{\Lambda_1} + \bar{\xi}_{\Lambda_1} \bar{\zeta}_{\Lambda_2} - \bar{\zeta}_{\Lambda_1} \bar{\xi}_{\Lambda_2}.$$

The *corrected* average current is given by

$$(8.7) \qquad \hat{w}_{\Lambda_1, \Lambda_2}(m) = \bar{w}_{\Lambda_1, \Lambda_2} - [\gamma_{\Lambda_2} - \gamma_{\Lambda_1}],$$

where

$$(8.8) \qquad \begin{aligned} \gamma_\Lambda = {} & 2m(1-m)\hat{\lambda}_\Lambda + m e^{-\lambda(m)}(\bar{\zeta}_\Lambda - \langle \bar{\zeta}_\Lambda \rangle_{\hat{\lambda}_\Lambda}) \\ & - (1-m)e^{\lambda(m)}(\bar{\xi}_\Lambda - \langle \bar{\xi}_\Lambda \rangle_{\hat{\lambda}_\Lambda}). \end{aligned}$$

Note that $\hat{w}_{\Lambda_1, \Lambda_2}$ is a function which depends on the variables $\alpha_x$ and $\eta_x$ for $x \in \Lambda_1 \cup \Lambda_2$, as well as an additional variable $m \in [0, 1]$. The form of the correction is chosen so that there is a $C < \infty$ such that if $\Lambda_1$ and $\Lambda_2$ are cubes of side length $\ell$, for all $m \in [0, 1]$, all $\Omega \subset \mathbb{Z}^d$ containing $\Lambda_1 \cup \Lambda_2$,

$$(8.9) \qquad \mathbb{E}[E[\{\hat{w}_{\Lambda_1, \Lambda_2}(m)\}^2 | \bar{\eta}_\Omega = m]] \leq C\ell^{-2d}.$$

We will prove this in the next section as Lemma 9.1.

We now set up the renormalization procedure.

Fix a large positive integer $M$ and let $\Omega_M$ be the union of the cube $\Lambda_M^1 = \{1, \ldots, M\}^d$ and its translate $\Lambda_M^2 = \tau_{Me}\Lambda_M^1$.

Let $r_i$, $i = 1, \ldots, M$, be defined through

$$(8.10) \qquad \frac{1}{M} Av_{i=1}^M Av_{j=M+1}^{2M} \sum_{k=i}^{j-1} a_k = Av_{i=1}^{2M} r_i a_i$$

for all $a_i$, $i = 1, \ldots, M$. Although it is elementary to write down a closed form expression for the $r_i$, we will only ever use (8.10). Next define

$$(8.11) \qquad \rho_x = r_{x \cdot e}, \qquad x \in \Omega_M.$$

Note that $Av_{x \in \Omega_M} \rho_x = 1$ and the average $Av_{x \in \Omega_M} \rho_x a_x$ of $a_x$, $x \in \Omega_M$ has the special property that if $a_x = A_{x+e} - A_x$, then

$$(8.12) \qquad Av_{x \in \Omega_M} \rho_x a_x = \frac{1}{M}[Av_{y \in \Lambda_M^2} A_y - Av_{x \in \Lambda_M^1} A_x].$$

LEMMA 8.1. *Let $M = K/\ell$ be an integer. Recall the definition (3.29) long jump current*

$$W_\ell^e = \ell^{-1} \bar{w}_{\Lambda^\ell, \tau_{\ell e}\Lambda^\ell}.$$

*Let*

(8.13) $$V_{K,\ell} = Av_{x \in \Omega_{K/\ell}} \rho_x \tau_{\ell x} W_\ell^e - W_K^e.$$

*There is a $C < \infty$ such that, for all $m \in [0,1]$, and positive integers $M$ and $\ell$,*

(8.14) $$K^d \mathbb{E}[E[V_{K,\ell}(-L_{\Omega_K})^{-1} V_{K,\ell} | \bar{\eta}_{\Omega_K} = m]] \leq CM^{d+2}\ell^{-d}.$$

PROOF. From (8.12),

(8.15) $$\ell^{-1} Av_{x \in \Omega_{K/\ell}} \rho_x \tau_{\ell x}[\tau_{\ell e}\gamma_\Lambda - \gamma_\Lambda]$$
$$= K^{-1}[Av_{y \in \Lambda_{K/\ell}^2} \tau_{\ell y}\gamma_\Lambda - Av_{x \in \Lambda_{K/\ell}^1} \tau_{\ell x}\gamma_\Lambda],$$

so we can write $Av_{y \in \Omega_{K/\ell}} \rho_x \tau_{\ell x} W_\ell^e - W_K^e = A_{K,\ell} + B_{K,\ell}$, where

(8.16) $$A_{K,\ell} = Av_{x \in \Omega_{K/\ell}} \rho_x \tau_{\ell x} \ell^{-1} \hat{w}_{\Lambda_\ell, \tau_{\ell e}\Lambda_\ell},$$

(8.17) $$B_{K,\ell} = Av_{x \in \Lambda_{K/\ell}^1, y \in \Lambda_{K/\ell}^2} K^{-1} \hat{w}_{\tau_{\ell x}\Lambda_\ell, \tau_{\ell y}\Lambda_\ell}.$$

The $\hat{w}$ are defined in (8.7). By Jensen's inequality,

(8.18) $$\mathbb{E}[E[A_{K,\ell}^2 | \bar{\eta}_{\Omega_K} = m]] \leq Av_{x \in \Omega_{K/\ell}} \rho_x \ell^{-2} \tau_{\ell x} \mathbb{E}[E[\hat{w}_\ell^2 | \bar{\eta}_{\Omega_K}]] \leq C\ell^{-2-2d}.$$

In the same way,

(8.19) $$\mathbb{E}[E[B_{K,\ell}^2 | \bar{\eta}_{\Omega_K} = m]] \leq CK^{-2}\ell^{-2d}.$$

By spectral gap,

(8.20) $$E[V_K(-L_{\Omega_K})^{-1} V_K | \bar{\eta}_{\Omega_K}] \leq CK^2 E[V_K^2 | \bar{\eta}_{\Omega_K}].$$

We conclude that the left-hand side of (8.14) is bounded by $CK^d \ell^{-2d}$. Since $K > \ell$, this is again bounded by $CK^{d+2}\ell^{-2d-2}$, which is the same as the right-hand side of (8.14). □

Now let $K = \ell M^N$ for some $N$ and $\ell_n = \ell M^n$. $\Omega_{\ell_{n+1}}$ is readily seen to be a union $\bigcup_{x \in \Omega_M} \tau_{\ell_n x} \Omega_{\ell_n}$ of copies of $\Omega_{\ell_n}$. So a point $x \in \Omega_{\ell_{n+1}}$ can be represented as $(x_1, x_2)$, where $x_1 \in \Omega_{\ell_n}$ and $x_2 \in \Omega_M$. Continuing in this way, we have

$$\Omega_K \simeq \Omega_\ell \times \Omega_M^N$$

and $\ell \mathbb{Z}^d \cap \Omega_K \simeq \Omega_M^N$.

If $x \in \Omega_M^n$, we can write $x = (x_1, \ldots, x_n)$, where $x_i \in \Omega_M$ and if $a_x$ is a function on $\Omega_M^n$, we can define the average

$$\text{(8.21)} \quad Av_{x \in \Omega_M^n} a_x \rho_x^{(n)} = Av_{x_1 \in \Omega_M} \cdots Av_{x_n \in \Omega_M} a_{(x_1,\ldots,x_n)} \rho_{x_1} \cdots \rho_{x_n}.$$

The main result of this section is the following:

THEOREM 4. *Fix $M$ a large integer and $K = \ell M^N$. Let*

$$\text{(8.22)} \quad R_{K,\ell} = Av_{x \in \ell \mathbb{Z}^d \cap \Omega_K} \tau_x W_\ell^e \rho_x^{(N)} - W_K^e.$$

*Then*

$$\text{(8.23)} \quad \lim_{\ell \to \infty} \lim_{N \to \infty} K^d \mathbb{E}[E[R_{K,\ell}, (-L_{\Omega_K})^{-1} R_{K,\ell} | \bar{\eta}_{\Omega_K}]] = 0.$$

PROOF. Let

$$\text{(8.24)} \quad R^{(n)} = Av_{x \in \Omega_M} \rho_x \tau_{\ell_n x} W_{\ell_n}^e - W_{\ell_{n+1}}^e$$

and

$$\text{(8.25)} \quad \bar{R}^{(n)} = Av_{x \in \Omega_M^{N-n-1}} \rho_x^{(N-n-1)} \tau_{\ell_{n+1} x} R^{(n)}$$

so that

$$\text{(8.26)} \quad R_{K,\ell} = \sum_{n=0}^{N-1} \bar{R}^{(n)}.$$

By the triangle inequality,

$$\text{(8.27)} \quad \begin{aligned} &(K^d \mathbb{E}[E[R_{K,\ell}(-L_{\Omega_K})^{-1} R_{K,\ell}]])^{1/2} \\ &\leq \sum_{n=0}^{N-1} (K^d \mathbb{E}[E[\bar{R}^{(n)}(-L_{\Omega_K})^{-1} \bar{R}^{(n)}]])^{1/2}. \end{aligned}$$

From the variational formula,

$$K^d E[\bar{R}^{(n)}(-L_{\Omega_K})^{-1} \bar{R}^{(n)} | \bar{\eta}_{\Omega_K}] \leq Av_{x \in \Omega_M^{N-n-1}} \rho_x^{(N-n-1)} \tau_{\ell_{n+1} x} \ell_{n+1}^d$$
$$\times \sup_m E[R^{(n)}(-L_{\Omega_{\ell_{n+1}}})^{-1} R^{(n)} | \bar{\eta}_{\Omega_{\ell_{n+1}}} = m].$$

By the previous lemma,

$$\text{(8.28)} \quad \ell_{n+1}^d \mathbb{E}[E[R^{(n)}(-L_{\Omega_{\ell_{n+1}}})^{-1} R^{(n)} | \bar{\eta}_{\Omega_{\ell_{n+1}}}]] \leq C M^{d+2} \ell_n^{-d}.$$

Hence, we have

$$\text{(8.29)} \quad \begin{aligned} (K^d \mathbb{E}[E[R_{K,\ell}, (-L_{\Omega_K}^{-1})^{-1} R_{K,\ell} | \bar{\eta}_{\Omega_{\ell_{n+1}}}]])^{1/2} &\leq C M^{(d+2)/2} \sum_{n=0}^{\infty} \ell_n^{-d/2} \\ &= C(M) \ell^{-d/2}, \end{aligned}$$

with $C(M) < \infty$. Let $\ell \to \infty$ to complete the proof. $\square$

Finally we apply the result to our particular problem. We need to show that the difference between (8.1) and

$$\int_0^T \varepsilon^d \sum_{x \in \mathbb{Z}^d/\varepsilon^{-1}\mathbb{Z}^d, e} J(\varepsilon x) \tau_x \varepsilon^{-1} \nu_{e_0, e}(\bar{\eta}_{\Lambda_K}) W_K^e \, dt \tag{8.30}$$

goes to 0 in $P_\varepsilon$-probability as $\varepsilon \to 0$ followed by $\ell \to \infty$, where $K = c_0 \varepsilon^{-2/(d+2)}$. First of all we would like to replace (8.1) by a term corresponding to the first term in (8.22). The difference is

$$\int_0^T \varepsilon^d \sum_{x \in \mathbb{Z}^d/\varepsilon^{-1}\mathbb{Z}^d, e} J(\varepsilon x) \tau_x \varepsilon^{-1} \nu_{e_0, e}(\bar{\eta}_{\Lambda_K})$$
$$\times [Av_{y \in \Lambda_K'} \tau_y W_\ell^e - Av_{x \in \ell \mathbb{Z}^d \cap \Omega_K} \tau_x W_\ell^e \rho_x^{(N)}] \, dt. \tag{8.31}$$

Performing a summation by parts, this can be rewritten as $\int_0^T \Gamma \, dt$, where

$$\Gamma = \varepsilon^d \sum_{x \in \mathbb{Z}^d/\varepsilon^{-1}\mathbb{Z}^d, e} Q_{K,\ell}(x) \varepsilon^{-1} \tau_x W_\ell^e \tag{8.32}$$

and

$$Q_{K,\ell}(x) = Av_{y \in \Lambda_K'} J(\varepsilon(x+y)) \nu_{e_0, e}(\tau_{x+y} \bar{\eta}_{\Lambda_K})$$
$$- Av_{y \in \ell \mathbb{Z}^d \cap \Omega_K} \rho_y^{(N)} J(\varepsilon(x+y)) \nu_{e_0, e}(\tau_{x+y} \bar{\eta}_{\Lambda_K}). \tag{8.33}$$

Now because $\nabla_{xy} \bar{\eta}_{\Lambda_K} = 0$ for $x, y \in \Lambda_K$,

$$\langle Q_{K,\ell}(x) \tau_x W_\ell^e f \rangle = \ell^{-1} Av_{y \in \tau_x \Lambda^\ell, z \in \tau_{x+e\ell} \Lambda^\ell} \langle Q_{K,\ell}(x)(\eta_z - \eta_y) \nabla_{yz} f \rangle. \tag{8.34}$$

Furthermore,

$$\langle Q_{K,\ell}(x)(\eta_z - \eta_y) \nabla_{yz} f \rangle = 2 \langle Q_{K,\ell}(x)(\eta_z - \eta_y) \sqrt{f} \nabla_{yz} \sqrt{f} \rangle$$
$$\leq 2 \langle Q_{K,\ell}^2(x) f \rangle^{1/2} \langle (\nabla_{yz} \sqrt{f})^2 \rangle^{1/2} \tag{8.35}$$

so

$$\langle \Gamma f \rangle \leq \varepsilon^{-1} Av_{x \in \mathbb{Z}^d/\varepsilon^{-1}\mathbb{Z}^d, e} \langle Q_{K,\ell}^2(x) f \rangle^{1/2}$$
$$\times \langle \ell^{-2} Av_{y \in \tau_x \Lambda^\ell, z \in \tau_{x+e\ell} \Lambda^\ell} (\nabla_{yz} \sqrt{f})^2 \rangle^{1/2}. \tag{8.36}$$

By the moving particle Lemma 5.2, there is a $C < \infty$ such that

$$\ell^{-2} Av_{y \in \tau_x \Lambda^\ell, z \in \tau_{x+e\ell} \Lambda^\ell} \langle (\nabla_{yz} \sqrt{f})^2 \rangle \leq C Av_{b \in \tau_x \Lambda^\ell \cup \tau_{x+e\ell} \Lambda^\ell} \langle (\nabla_b \sqrt{f})^2 \rangle, \tag{8.37}$$

where the average on the right-hand side is over nearest neighbor bonds only. Hence,

$$
\begin{aligned}
\langle \Gamma f \rangle &- \varepsilon^{d-2} D(\sqrt{f}) \\
&\leq C A v_{x,e} \varepsilon^{-1} \langle Q_{K,\ell}^2(x) f \rangle^{1/2} \langle A v_{|y-x|\leq 2\ell}(\nabla_{x,x+e}\sqrt{f})^2 \rangle^{1/2} \\
&\quad - \varepsilon^{-2} \langle (\nabla_{x,x+e}\sqrt{f})^2 \rangle \\
&\leq C' A v_{x \in \mathbb{Z}^d/\varepsilon^{-1}\mathbb{Z}^d, e} \langle Q_{K,\ell}^2(x) f \rangle - \tfrac{1}{2} \varepsilon^{d-2} D(\sqrt{f}),
\end{aligned}
\tag{8.38}
$$

which vanishes in the limit of small $\varepsilon$ by the two block estimate (5.6).

Summarizing the results so far, we have shown that the difference between (3.8) and (3.15) vanishes in the limit $\varepsilon \downarrow 0$, assuming (8.9). This is proved in the next section. Following that, we still have to show that the difference of (3.15) and (3.10) is small to obtain the hydrodynamic equation.

**9. Variance estimate for the corrected average current.** In this section we give the proof of (8.9) which is the key to the renormalization procedure of the previous section.

LEMMA 9.1.  *Let $\hat{w}_{\Lambda_1,\Lambda_2}(m)$ be as in (8.7). There is a constant $C < \infty$ such that, for all $m \in [0,1]$ and positive integers $K$ and $\ell$ with $K > 2\ell$,*

$$
\mathbb{E}[E[\{\hat{w}_{\Lambda_1,\Lambda_2}(m)\}^2 | \bar{\eta}_\Omega = m]] \leq C \ell^{-2d}. \tag{9.1}
$$

PROOF. Note that $\langle \xi_x \rangle_\lambda = \frac{1}{1+e^{\lambda+\alpha_x}} e^{\alpha_x} = e^{-\lambda} \langle \eta_x \rangle_\lambda$ and $\langle \zeta_x \rangle_\lambda = \frac{e^{\lambda+\alpha_x}}{1+e^{\lambda+\alpha_x}} \times e^{-\alpha_x} = e^\lambda \langle 1-\eta_x \rangle_\lambda$, so that

$$
\langle \bar{\xi}_\Lambda \rangle_{\hat{\lambda}_\Lambda} = e^{-\hat{\lambda}_\Lambda} \bar{\eta}_\Lambda, \qquad \langle \bar{\zeta}_\Lambda \rangle_{\hat{\lambda}_\Lambda} = e^{\hat{\lambda}_\Lambda} (1-\bar{\eta}_\Lambda). \tag{9.2}
$$

This gives

$$
\bar{w}_{\Lambda_1,\Lambda_2} - 2m(1-m)(\hat{\lambda}_{\Lambda_2} - \hat{\lambda}_{\Lambda_1}) = B_1 + B_2 - B_3, \tag{9.3}
$$

where

$$
B_1 = \bar{\xi}_{\Lambda_1} \bar{\zeta}_{\Lambda_2} - \langle \bar{\xi}_{\Lambda_1} \rangle_{\hat{\lambda}_{\Lambda_1}} \langle \bar{\zeta}_{\Lambda_2} \rangle_{\hat{\lambda}_{\Lambda_2}} - \bar{\zeta}_{\Lambda_1} \bar{\xi}_{\Lambda_2} + \langle \bar{\zeta}_{\Lambda_1} \rangle_{\hat{\lambda}_{\Lambda_1}} \langle \bar{\xi}_{\Lambda_2} \rangle_{\hat{\lambda}_{\Lambda_2}}, \tag{9.4}
$$

$$
B_2 = (2m(1-m) - \bar{\eta}_{\Lambda_1}(1-\bar{\eta}_{\Lambda_2}) - \bar{\eta}_{\Lambda_2}(1-\bar{\eta}_{\Lambda_1}))(\hat{\lambda}_{\Lambda_2} - \hat{\lambda}_{\Lambda_1}), \tag{9.5}
$$

and, setting $\phi(x) = e^x - 1 - x$,

$$
B_3 = \bar{\eta}_{\Lambda_1}(1-\bar{\eta}_{\Lambda_2}) \phi(\hat{\lambda}_{\Lambda_2} - \hat{\lambda}_{\Lambda_1}) - \bar{\eta}_{\Lambda_2}(1-\bar{\eta}_{\Lambda_1}) \phi(-(\hat{\lambda}_{\Lambda_2} - \hat{\lambda}_{\Lambda_1})). \tag{9.6}
$$

Decomposing $B_1 = B_{1,1} + B_{1,2} + B_{1,3}$,

$$B_{1,1} = (\bar\xi_{\Lambda_1} - \langle\bar\xi_{\Lambda_1}\rangle_{\hat\lambda_{\Lambda_1}})(\bar\zeta_{\Lambda_2} - \langle\bar\zeta_{\Lambda_2}\rangle_{\hat\lambda_{\Lambda_2}})$$
$$- (\bar\xi_{\Lambda_2} - \langle\bar\xi_{\Lambda_2}\rangle_{\hat\lambda_{\Lambda_2}})(\bar\zeta_{\Lambda_1} - \langle\bar\zeta_{\Lambda_1}\rangle_{\hat\lambda_{\Lambda_1}}),$$

$$B_{1,2} = (\langle\bar\xi_{\Lambda_1}\rangle_{\hat\lambda_{\Lambda_1}} - me^{-\lambda(m)})(\bar\zeta_{\Lambda_2} - \langle\bar\zeta_{\Lambda_2}\rangle_{\hat\lambda_{\Lambda_2}})$$
$$- (\langle\bar\xi_{\Lambda_2}\rangle_{\hat\lambda_{\Lambda_2}} - me^{-\lambda(m)})(\bar\zeta_{\Lambda_1} - \langle\bar\zeta_{\Lambda_1}\rangle_{\hat\lambda_{\Lambda_1}})$$

(9.7)
$$+ (\langle\bar\zeta_{\Lambda_2}\rangle_{\hat\lambda_{\Lambda_2}} - (1-m)e^{\lambda(m)})(\bar\xi_{\Lambda_1} - \langle\bar\xi_{\Lambda_1}\rangle_{\hat\lambda_{\Lambda_1}})$$
$$- (\langle\bar\zeta_{\Lambda_1}\rangle_{\hat\lambda_{\Lambda_1}} - (1-m)e^{\lambda(m)})(\bar\xi_{\Lambda_2} - \langle\bar\xi_{\Lambda_2}\rangle_{\hat\lambda_{\Lambda_2}}),$$

$$B_{1,3} = me^{-\lambda(m)}[(\bar\zeta_{\Lambda_2} - \langle\bar\zeta_{\Lambda_2}\rangle_{\hat\lambda_{\Lambda_2}}) - (\bar\zeta_{\Lambda_1} - \langle\bar\zeta_{\Lambda_1}\rangle_{\hat\lambda_{\Lambda_1}})]$$
$$+ (1-m)e^{\lambda(m)}[(\bar\xi_{\Lambda_1} - \langle\bar\xi_{\Lambda_1}\rangle_{\hat\lambda_{\Lambda_1}}) - (\bar\xi_{\Lambda_2} - \langle\bar\xi_{\Lambda_2}\rangle_{\hat\lambda_{\Lambda_2}})].$$

We will obtain a bound $\mathbb{E}[E[B^2|\bar\eta_\Omega]] \leq C\ell^{-2d}$ for each of $B_{1,1}, B_{1,2}, B_2$ and $B_3$. Note that $B_{1,3}$ is the extra term appearing in $\gamma_{\Lambda_2} - \gamma_{\Lambda_1}$. Through the proof, $C$ will stand for a finite constant independent of $\bar\eta_\Omega$, though its meaning will change from line to line.

We start with $B_{1,1}$. By the Schwarz inequality,

$$E[(\bar\xi_{\Lambda_1} - E[\bar\xi_{\Lambda_1}|\bar\eta_{\Lambda_1}])^2(\bar\zeta_{\Lambda_2} - E[\bar\zeta_{\Lambda_2}|\bar\eta_{\Lambda_2}])^2|\bar\eta_\Omega] \leq C\ell^{-2d}$$

and the same for the analogue of the second term of $B_{1,1}$. By the equivalence of ensembles (see Appendix 2 of [11]),

$$|\langle\bar\xi_{\Lambda_1}\rangle_{\hat\lambda_{\Lambda_1}} - E[\bar\xi_{\Lambda_1}|\bar\eta_{\Lambda_1}]| \leq C\ell^{-d}.$$

Hence, it is not hard to compute

$$E[B_{1,1}^2|\bar\eta_\Omega] \leq C\ell^{-2d}.$$

When we try to do the same thing for $E[B_{1,2}^2|\bar\eta_\Omega]$, we get terms like $E[(\bar\zeta_{\Lambda_2} - E[\bar\zeta_{\Lambda_2}|\bar\eta_{\Lambda_2}])^4|\bar\eta_\Omega]$, which are bounded by $C\ell^{-2d}$, but also a term $E[(\langle\bar\xi_{\Lambda_1}\rangle_{\hat\lambda_{\Lambda_1}} - me^{-\lambda(m)})^4|\bar\eta_\Omega]$ for which the same argument does not work. Recall $\langle\bar\xi_{\Lambda_1}\rangle_{\hat\lambda_{\Lambda_1}} = \bar\eta_{\Lambda_1}e^{-\hat\lambda_{\Lambda_1}}$ and rewrite the difference as two terms

$$[\bar\eta_{\Lambda_1}e^{-\lambda(\bar\eta_{\Lambda_1})} - me^{-\lambda(m)}] + [\bar\eta_{\Lambda_1}e^{-\lambda(\bar\eta_{\Lambda_1})}(e^{-(\hat\lambda_{\Lambda_1}(\bar\eta_{\Lambda_1})-\lambda(\bar\eta_{\Lambda_1}))} - 1)].$$

The function $\rho e^{-\lambda(\rho)}$ is bounded and Lipschitz and $\mathbb{E}[\mathbb{E}[(\bar\eta_{\Lambda_1} - m)^4|\bar\eta_\Omega = m]] \leq C\ell^{-2d}$. Hence, for the first term, we have the required bound $\mathbb{E}[\mathbb{E}[(\bar\eta_{\Lambda_1}e^{-\lambda(\bar\eta_{\Lambda_1})} - me^{-\lambda(m)})^4|\bar\eta_\Omega = m]] \leq C\ell^{-2d}$. For the second term, note that $\bar\eta_{\Lambda_1}e^{-\lambda(\bar\eta_{\Lambda_1})}$ is uniformly bounded. For $|x| \leq A$, $|e^x - 1| \leq C|x|$ so

$\mathbb{E}[E[\mathbb{1}(|\hat{\lambda}_{\Lambda_1}(\bar{\eta}_{\Lambda_1}) - \lambda(\bar{\eta}_{\Lambda_1})| \leq A)(e^{-(\hat{\lambda}_{\Lambda_1}(\bar{\eta}_{\Lambda_1}) - \lambda(\bar{\eta}_{\Lambda_1}))} - 1)^4 |\bar{\eta}_\Omega]$ is bounded above by a constant multiple of $\mathbb{E}[E[(\hat{\lambda}_{\Lambda_i}(\bar{\eta}_{\Lambda_i}) - \lambda(\bar{\eta}_{\Lambda_i}))^4 |\bar{\eta}_\Omega]]$. By Lemma 10.3,

$$\mathbb{E}[(\hat{\lambda}_{\Lambda_i}(\bar{\eta}_{\Lambda_i}) - \lambda(\bar{\eta}_{\Lambda_i}))^4] \leq C\ell^{-2d} \tag{9.8}$$

for $i = 1, 2$. On the other hand, $\bar{\eta}_{\Lambda_1} e^{-(\hat{\lambda}_{\Lambda_1}(\bar{\eta}_{\Lambda_1}))}$ and $\bar{\eta}_{\Lambda_1} e^{-\lambda(\bar{\eta}_{\Lambda_1})}$ are both uniformly bounded and by Chebyshev's inequality and (9.8), $\mathbb{E}[E[\mathbb{1}(|\hat{\lambda}_{\Lambda_1}(\bar{\eta}_{\Lambda_1}) - \lambda(\bar{\eta}_{\Lambda_1})| > A)|\bar{\eta}_\Omega] \leq CA^{-4}\ell^{-2d}$. We conclude that

$$\mathbb{E}[E[B_{1,2}^2|\bar{\eta}_\Omega]] \leq C\ell^{-2d}.$$

Turning to $B_2$, it is not hard to see that

$$\mathbb{E}[E[((2m(1-m) - \bar{\eta}_{\Lambda_1}(1 - \bar{\eta}_{\Lambda_2}) - \bar{\eta}_{\Lambda_2}(1 - \bar{\eta}_{\Lambda_1})))^4|\bar{\eta}_\Omega]] \leq C\ell^{-2d}.$$

We claim that

$$\mathbb{E}[E[(\hat{\lambda}_{\Lambda_2} - \hat{\lambda}_{\Lambda_1})^4|\bar{\eta}_\Omega]] \leq C\ell^{-2d}. \tag{9.9}$$

Then by Schwarz's inequality, we conclude that $\mathbb{E}[E[B_2^2|\bar{\eta}_\Omega]] \leq C\ell^{-2d}$. To prove (9.9), rewrite $\hat{\lambda}_{\Lambda_2} - \hat{\lambda}_{\Lambda_1} = \hat{\lambda}_{\Lambda_2}(\bar{\eta}_{\Lambda_2}) - \hat{\lambda}_{\Lambda_1}(\bar{\eta}_{\Lambda_1})$ as

$$[\hat{\lambda}_{\Lambda_2}(\bar{\eta}_{\Lambda_2}) - \lambda(\bar{\eta}_{\Lambda_2})] + [\lambda(\bar{\eta}_{\Lambda_2}) - \lambda(\bar{\eta}_{\Lambda_1})] - [\hat{\lambda}_{\Lambda_1}(\bar{\eta}_{\Lambda_1}) - \lambda(\bar{\eta}_{\Lambda_1})].$$

The first and third terms can be handled by (9.8). Using the fact that the variance of $\eta_x$ is bounded by a constant times $1/\lambda'$, it is not hard to check that also $\mathbb{E}[E[(\lambda(\bar{\eta}_{\Lambda_2}) - \lambda(\bar{\eta}_{\Lambda_1}))^4|\bar{\eta}_\Omega]] \leq C\ell^{-2d}$, proving (9.9).

Finally we consider $B_3$. Let $X = \hat{\lambda}_{\Lambda_2} - \hat{\lambda}_{\Lambda_1}$. From (9.9), we know that $\mathbb{E}[E[X^4|\bar{\eta}_\Omega]] \leq C\ell^{-2d}$, so it will suffice to bound $\mathbb{E}[E[B_3^2|\bar{\eta}_\Omega]] \leq C\mathbb{E}[E[X^4|\bar{\eta}_\Omega]]$. There is a constant $C = C(A) < \infty$ so that if $|x| \leq A$, then $|\phi(x)| \leq Cx^2$. So $\mathbb{1}(|X| \leq A)\phi^2(X) \leq CX^4$ and, hence,

$$\mathbb{E}[E[\mathbb{1}(|X| \leq A)B_3^2|\bar{\eta}_\Omega]] \leq C\mathbb{E}[E[X^4|\bar{\eta}_\Omega]].$$

On the other hand, it is not hard to check that

$$Y = \bar{\eta}_{\Lambda_1}(1 - \bar{\eta}_{\Lambda_2})(e^{\hat{\lambda}_{\Lambda_2} - \hat{\lambda}_{\Lambda_1}} - 1) - \bar{\eta}_{\Lambda_2}(1 - \bar{\eta}_{\Lambda_1})(e^{\hat{\lambda}_{\Lambda_1} - \hat{\lambda}_{\Lambda_2}} - 1)$$

is uniformly bounded. Hence,

$$\mathbb{E}[E[\mathbb{1}(|X| \geq A)Y^2|\bar{\eta}_\Omega]] \leq C\mathbb{E}[E[\mathbb{1}(|X| \geq A)|\bar{\eta}_\Omega]]$$

and by Chebyshev's inequality, the last term is bounded above by $A^{-4}\mathbb{E}[E[X^4|\bar{\eta}_\Omega]]$. Finally, $\mathbb{E}[E[\mathbb{1}(|X| \geq A)X^2|\bar{\eta}_\Omega]] \leq 2\mathbb{E}[E[\mathbb{1}(|X| \geq A)|\bar{\eta}_\Omega]] + 2\mathbb{E}[E[X^4|\bar{\eta}_\Omega]]$. The first part is bounded by $2A^{-4}\mathbb{E}[E[X^4|\bar{\eta}_\Omega]]$ by Chebyshev's inequality again. $\square$

LEMMA 9.2.

$$|\langle\bar{\xi}_{\Lambda_1}\rangle_{\hat{\lambda}_{\Lambda_1}} - E[\bar{\xi}_{\Lambda_1}|\bar{\eta}_{\Lambda_1}]| \leq C\ell^{-d}.$$

PROOF. Fix $\lambda$ and let $P$ denote the product measure with $P(\eta_x = 1) = \frac{e^{\alpha_x + \lambda}}{1 + e^{\alpha_x + \lambda}}$ and $E$ the corresponding expectation. We claim first that

$$\text{(9.10)} \qquad \frac{P(\sum_{x \in \Lambda} \eta_x = M - 1)}{P(\sum_{x \in \Lambda} \eta_x = M)} = \frac{|\Lambda|}{|\Lambda| - M + 1} e^{-\lambda} E\left[\bar{\zeta}_\Lambda \Big| \sum_{x \in \Lambda} \eta_x = M\right].$$

Now

$$\langle \bar{\zeta}_\Lambda \rangle_{\hat{\lambda}_\Lambda} - E\left[\bar{\zeta}_\Lambda \Big| \sum_{x \in \Lambda} \eta_x = M\right]$$

$$= \sum_{x \in \Lambda} \frac{e^{\hat{\lambda}_\Lambda}}{1 + e^{\alpha_x + \hat{\lambda}_\Lambda}} \left[1 - \frac{P(\sum_{y \in \Lambda - \{x\}} \eta_y = M - 1)}{P(\sum_{y \in \Lambda} \eta_y = M)}\right]$$

$$= \sum_{x \in \Lambda} \frac{e^{\hat{\lambda}_\Lambda}}{1 + e^{\alpha_x + \hat{\lambda}_\Lambda}}$$

$$\times \left[1 - \left\{P(\eta_x = 1) + P(\eta_x = 0) \frac{P(\sum_{y \in \Lambda - \{x\}} \eta_y = M)}{P(\sum_{y \in \Lambda - \{x\}} \eta_y = M - 1)}\right\}^{-1}\right].$$

Now

$$\text{(9.11)} \qquad \frac{P(\sum_{y \in \Lambda} \eta_y = M)}{P(\sum_{y \in \Lambda - \{x\}} \eta_y = M - 1)}$$

$$= \frac{e^{\alpha_x + \hat{\lambda}_\Lambda}}{1 + e^{\alpha_x + \hat{\lambda}_\Lambda}} + \frac{1}{1 + e^{\alpha_x + \hat{\lambda}_\Lambda}} \frac{P(\sum_{y \in \Lambda - \{x\}} \eta_y = M)}{P(\sum_{y \in \Lambda - \{x\}} \eta_y = M - 1)}.$$

By (9.10),

$$\text{(9.12)} \qquad \frac{P(\sum_{y \in \Lambda - \{x\}} \eta_y = M)}{P(\sum_{y \in \Lambda - \{x\}} \eta_y = M - 1)}$$

$$= \left\{\frac{|\Lambda| - 1}{|\Lambda| - M} e^{-\hat{\lambda}_\Lambda} E\left[\bar{\zeta}_{\Lambda - \{x\}} \Big| \sum_{y \in \Lambda - \{x\}} \eta_y = M\right]\right\}^{-1}.$$

By the inductive hypothesis,

$$\text{(9.13)} \qquad E\left[\bar{\zeta}_{\Lambda - \{x\}} \Big| \sum_{y \in \Lambda - \{x\}} \eta_y = M\right] = \langle \bar{\zeta}_{\Lambda - \{x\}} \rangle_{\hat{\lambda}_{\Lambda - \{x\}}} + \frac{C}{|\Lambda| - 1}$$

$$= e^{\hat{\lambda}_{\Lambda - \{x\}}} \frac{|\Lambda| - M}{|\Lambda| - 1} + \frac{C}{|\Lambda| - 1}. \qquad \square$$

**10. Fick's law.** Our goal in this section is to prove that, for any smooth $J$,

$$\int_0^T \varepsilon^d \sum_{x \in \varepsilon^{-1}\mathbb{Z}^d/\mathbb{Z}^d} J(\varepsilon x)[\varepsilon^{-1}\tau_x \Theta_K^e - D_{e0e}(\bar{\eta}_x^{\delta_1 \varepsilon^{-1}})(2\delta_2)^{-1}$$
$$\times (\bar{\eta}_{x+\delta_2\varepsilon^{-1}e}^{\delta_1\varepsilon^{-1}} - \bar{\eta}_{x-\delta_2\varepsilon^{-1}e}^{\delta_1\varepsilon^{-1}})] \, dt \quad (10.1)$$

vanishes in $P_\varepsilon$ probability, as $\delta_1, \delta_2 \to 0$. $\Theta_K^e$ is given is (3.28).

Let $m_1$ and $m_2$ be the particle densities on $\Lambda_K^1 = \Lambda_K$ and $\Lambda_K^2 = \tau_{Ke}\Lambda_K$. Let $\hat{\lambda}_i(m) = \hat{\lambda}_\Lambda^i(m)$, $i = 1, 2$, be the empirical chemical potential on $\Lambda_K^i$ as in (8.4). Let

$$\hat{F}(m_1, m_2) = \langle W_K^e \rangle_{\hat{\lambda}_1(m_1), \hat{\lambda}_2(m_2)}. \quad (10.2)$$

The expectation is with respect to the product measure on $\Lambda_1$ with chemical potential $\lambda_1(m_1)$ and on $\Lambda_2$ with chemical potential $\Lambda_2$:

$$Z^{-1} \exp\left\{\sum_{x \in \Lambda_1} (\alpha_x + \lambda_1(m_1))\eta_x + \sum_{y \in \Lambda_2} (\alpha_y + \lambda_2(m_2))\eta_y\right\}. \quad (10.3)$$

The function $\hat{F}(m_1, m_2)$, which depends on the densities $m_1$ and $m_2$, as well as on the field configurations on the blocks $\Lambda_K^1$ and $\Lambda_K^2$, is given explicitly by

$$\hat{F}(m_1, m_2) = K^{-1}[(e^{\hat{\lambda}_2(m_2) - \hat{\lambda}_1(m_1)} - 1)m_1(1 - m_2)$$
$$- (e^{\hat{\lambda}_1(m_1) - \hat{\lambda}_2(m_2)} - 1)m_2(1 - m_1)]. \quad (10.4)$$

Consider as well the variant of $\hat{F}$ where we use the annealed chemical potential instead of the empirical chemical potential

$$F(m_1, m_2) = K^{-1}[(e^{\lambda(m_2) - \lambda(m_1)} - 1)m_1(1 - m_2)$$
$$- (e^{\lambda(m_1) - \lambda(m_2)} - 1)m_2(1 - m_1)]. \quad (10.5)$$

For $m_1 \neq m_2$, consider the quotient $\Phi(m_1, m_2) = F(m_1, m_2)/(K^{-1}(m_2 - m_1))$. From the boundedness of the field $\alpha$, it is bounded above and below away from 0 and Lipschitz. When $m_1 = m_2 = m$, we have $\Phi(m, m) = 2\lambda'(m)m(1 - m)$. We call $G(m_1, m_2) = 1/\Phi(m_1, m_2)$. $G$ defined in this way is also bounded and uniformly Lipschitz on $[0, 1] \times [0, 1]$.

We will prove (10.1) is several steps. The first thing we want to do is replace the term

$$Av_x J(\varepsilon x)\varepsilon^{-1}\tau_x \Theta_K^e = Av_x \varepsilon^{-1} J(\varepsilon x)\tau_x \nu(\bar{\eta}_{\Lambda_K})W_K^e \quad (10.6)$$

in (10.1) by

$$(10.7) \quad Av_x \varepsilon^{-1} \tau_x \Gamma W_K^e, \qquad \Gamma = [Av_y \rho_y^{\delta_1,\delta_2} J(\varepsilon y) \tau_y D_{e_0 e}(\bar{\eta}_0^{\delta_1 \varepsilon^{-1}})]G,$$

where the weights $\rho_y^{\delta_1,\delta_2}$ are defined on the convex hull of two boxes of side length $\varepsilon^{-1}\delta_1$ whose centers are separated by the vector $\varepsilon^{-1}\delta_2 e$ and $\sum_y \rho_y^{\delta_1,\delta_2} = 1$. Here $\nu(m) = \sigma_{e_0 e}(m)/m(1-m)$, $G = G(m_1, m_2)$ and we have abused our definitions mildly by writing $\tau_x J(\varepsilon y) = J(\varepsilon(x+y))$.

The weights $\rho_y^{\delta_1,\delta_2}$ are defined as follows. We can assume without loss of generality that $K$ divides both $\varepsilon^{-1}\delta_1$ and $\varepsilon^{-1}\delta_2$ evenly. Divide the box centered at the origin of side length $\varepsilon^{-1}\delta_1$ into boxes of side length $K$ and label their centers $\beta$. For any $a_y$, $y \in \mathbb{Z}^d$, let

$$(10.8) \quad Av_y \rho_y^{\delta_1,\delta_2} a_y = Av_\beta (2\delta_2)^{-1} \varepsilon K \sum_{i=-\delta_2 \varepsilon^{-1} K^{-1}}^{\delta_2 \varepsilon^{-1} K^{-1}-1} a_{\beta+(i-1/2)Ke}.$$

This average has the property that if $a_y = \bar{\eta}_{y+Ke/2}^K - \bar{\eta}_{y-Ke/2}^K$,

$$(\varepsilon K)^{-1} Av_y \rho_y^{\delta_1,\delta_2} a_y = (2\delta_2)^{-1} (\bar{\eta}_{\delta_2 \varepsilon^{-1} e}^{\delta_1 \varepsilon^{-1}} - \bar{\eta}_{-\delta_2 \varepsilon^{-1} e}^{\delta_1 \varepsilon^{-1}}).$$

A summation by parts gives

$$\varepsilon^{d-1} K^{-1} \sum_x \tau_x \Gamma F = \varepsilon^d \sum_x J(\varepsilon x) \tau_x D_{e_0 e}(\bar{\eta}_0^{\delta_1 \varepsilon^{-1}}) \left( \frac{\bar{\eta}_{\delta_2 \varepsilon^{-1} e}^{\delta_1 \varepsilon^{-1}} - m_{-\delta_2 \varepsilon^{-1} e}^{\delta_1 \varepsilon^{-1}}}{2\delta_2} \right),$$

which is the right-hand side of (10.1). In other words, if we can replace $\nu$ by $\Gamma$ as described in (10.6) to (10.7), and then replace $W_K$ by $K^{-1}F$, we end up with the right-hand side of (10.1).

We start with (10.6) to (10.7). The definition of constants $C$ will change from line to line, but will always denote a finite constant independent of the parameters $\varepsilon$, $\delta_1$ and $\delta_2$. Let $\Psi = J\nu(\eta_0^{2K}) - \Gamma$. We have

$$(10.9) \quad \langle W_K \Psi f \rangle = K^{-1} Av_{x \in \Lambda_K^1, y \in \Lambda_K^2} \langle (\nabla_{xy} \eta)(\nabla_{xy} \Psi f) \rangle.$$

Now $\nabla_{xy} \Psi f = \Psi \nabla_{xy} f + f(T_{xy}\eta) \nabla_{xy} \Psi$. The first piece $\Psi \nabla_{xy} f$ is rather standard. Write $\nabla_{xy} f = \nabla_{xy}^+ \sqrt{f} \nabla_{xy} \sqrt{f}$, where $\nabla_{xy}^+ f = f(T_{xy}\eta) + f(\eta)$. Changing variables $T_{xy}\eta \mapsto \eta$, we can write $\langle \Psi \nabla_{xy} \eta \nabla_{xy} f \rangle = \langle (\nabla_{xy}^+ \Psi) \sqrt{f} \nabla_{xy} \eta \nabla_{xy} \sqrt{f} \rangle$. By Schwarz's inequality and the moving particle lemma, for any $q > 0$,

$$(10.10) \quad Av_{x \in \Lambda_K^1, y \in \Lambda_K^2} \langle \Psi \nabla_{xy} \eta \nabla_{xy} f \rangle \leq Cq \langle \Psi^2 f \rangle + K^2 q^{-1} \sum_{x, x+e \in \Lambda_K} \langle (\nabla_{x, x+e} \sqrt{f})^2 \rangle.$$

In the second piece, $f(T_{xy}\eta)\nabla_{xy}\Psi$,

$$\begin{aligned}(10.11)\quad \nabla_{xy}\Psi &= G(T_{xy}\eta)\nabla_{xy}Av_u\rho_u^{\delta_1,\delta_2}JD_{e_0e}(\bar\eta_u^{\delta_1\varepsilon^{-1}}) \\ &\quad + Av_u\rho_u^{\delta_1,\delta_2}JD_{e_0e}(\bar\eta_u^{\delta_1\varepsilon^{-1}})\nabla_{xy}G.\end{aligned}$$

From the definition (10.8) of the weights $\rho_u^{\delta_1,\delta_2}$ in the first term above, there are only a finite number of terms in the sum for which $\nabla_{xy}D_{e_0e}(\bar\eta_u^{\delta_1\varepsilon^{-1}})$ is nonzero and those that do give a term of the form $GJ(\delta_2)^{-1}\varepsilon K(D_{e_0e}(\bar\eta_u^{\delta_1\varepsilon^{-1}} \pm (\delta_1)^{-1}\varepsilon) - D_{e_0e}(\bar\eta_u^{\delta_1\varepsilon^{-1}}))$. For fixed $\delta_1,\delta_2 > 0$, this is $o(\varepsilon K)$. As for the second part, since $G$ is uniformly Lipschitz, $|\nabla_{xy}G| \leq CK^{-d}$. So we have

$$(10.12)\qquad Av_{x\in\Lambda_K^1, y\in\Lambda_K^2}\langle f\nabla_{xy}\eta\nabla_{xy}\Psi\rangle = o(\varepsilon K).$$

Hence, taking $q = \varepsilon K$ in (10.10),

$$\begin{aligned}\varepsilon^{-1}&Av_{x\in\mathbb{Z}^d/\varepsilon^{-1}\mathbb{Z}^d}\langle\tau_x\Psi W_K f\rangle - \varepsilon^{-2}Av_{x,e}\langle(\nabla_{xx+e}\sqrt f)^2\rangle \\ &\leq CAv_x\langle(\tau_x\Psi)^2 f\rangle - \tfrac{1}{2}\varepsilon^{-2}Av_{x,e}\langle(\nabla_{xx+e}\sqrt f)^2\rangle + o(1).\end{aligned}$$

The first terms on the right-hand side vanish as $\varepsilon \to 0$ followed by $\delta_1 \to 0$ and $\delta_2 \to 0$ by the two block estimate. To conclude, we have

$$(10.13)\quad \limsup_{\delta_2,\delta_1,\varepsilon}\sup_f\{\varepsilon^{-1}Av_x\langle\tau_x\Psi W_K f\rangle - \varepsilon^{-2}Av_{x,e}\langle(\nabla_{xx+e}\sqrt f)^2\rangle\} = 0,$$

which shows that we can replace (10.6) by (10.7).

The next step is to replace $W_K$ by $F(m_1,m_2)$. We do it in two steps: First replace $W_K$ by $\hat F$ (Lemma 10.2) and then $\hat F$ by $F$ (Lemma 10.4). Before proving Lemma 10.2, we need a preliminary estimate which shows that the integration by parts property of our long jump currents is preserved under conditional expectations.

LEMMA 10.1.  *Let $\Omega = A \cup B$ be subsets of $\mathbb{Z}^d$ with $\max_{x\in A, y\in B}|y-x| \leq K$. Let $x \in A$ and $y \in B$ and $w_{xy}$ be the current given by (3.1). Then there is a $C < \infty$ so that*

$$\begin{aligned}\langle E[w_{xy}|\bar\eta_A]f\rangle &\leq CK\sum_{u\in A}\langle(\nabla_{uy}\sqrt f)^2\rangle, \\ \langle E[w_{xy}|\bar\eta_B]f\rangle &\leq CK\sum_{v\in B}\langle(\nabla_{xv}\sqrt f)^2\rangle, \\ \langle E[w_{xy}|\bar\eta_A,\bar\eta_B]f\rangle &\leq CK\sum_{u\in A, v\in B}\langle(\nabla_{uv}\sqrt f)^2\rangle.\end{aligned}$$

PROOF. Integrating by parts,

$$(10.14)\qquad \langle E[w_{xy}|\bar\eta_A]f\rangle = \langle(\eta_y - \eta_x)\nabla_{xy}E[f|\bar\eta_A]\rangle.$$

There exist bounded $a_u$, $u \in A$ so that

$$\nabla_{xy} E[f|\bar{\eta}_A] = \sum_{u \in A} E[a_u \nabla_{uy} f|\bar{\eta}_A]. \tag{10.15}$$

In fact,

$$a_u = \eta_y(1-\eta_u)\frac{P(\eta^u|\sum_{z\in A}\eta_z + 1)}{P(\eta|\sum_{z\in A}\eta_z)} + (1-\eta_y)\eta_u\frac{P(\eta^u|\sum_{z\in A}\eta_z - 1)}{P(\eta|\sum_{z\in A}\eta_z)},$$

where $\eta^u$ is the configuration changed only at $u$. Hence,

$$\langle E[w_{xy}|\bar{\eta}_A]f\rangle = \sum_{u \in A} \langle b_u \nabla_{uy} f\rangle, \tag{10.16}$$

where $b_u = a_u E[(\eta_y - \eta_x)|\bar{\eta}_A]$. The result for $E[w_{xy}|\bar{\eta}_A]$ follows by Schwarz's inequality. The other results are proved in exactly the same way. □

LEMMA 10.2. *Let*

$$\Omega_1 = Av_x \tau_x \Gamma[W_K - \hat{F}(m_1, m_2)]. \tag{10.17}$$

*For any $\gamma > 0$,*

$$\limsup_{\delta_2 \to 0} \limsup_{\delta_1 \to 0} \limsup_{\varepsilon \to 0} \sup_{f} \{\varepsilon^{-1}\langle f\Omega_1\rangle \tag{10.18}$$
$$- \gamma\varepsilon^{-2} Av_{x,e}\langle(\nabla_{x,x+e}\sqrt{f})^2\rangle\} \leq 0.$$

PROOF. Let $F_1 = E[W_K|m_1]$ and $F_2 = E[W_K|m_2]$. By Lemma 10.1, $W_K$, $F_1$, $F_2$ and $F$ all satisfy (3.22). From these estimates and the two block estimate, we can replace $\Gamma$ by $\tilde{\Gamma}$, where $\tilde{\Gamma}\phi(\varepsilon y)D_{e_0 e}(\bar{\eta}^K))G(m_1, m_2)$. We will write $\Omega_1$ as $K^{-1} Av_x \tilde{\Gamma}\tau_x[(W_K - F_1 - F_2 + F) + (F_1 + F_2 - 2F)]$. Now

$$W_K - F_1 - F_2 + F = (\bar{\xi}_{\Lambda_1} - E[\bar{\xi}_{\Lambda_1}|\bar{\eta}_{\Lambda_1}])(\bar{\zeta}_{\Lambda_2} - E[\bar{\zeta}_{\Lambda_2}|\bar{\eta}_{\Lambda_2}])$$
$$- (\bar{\zeta}_{\Lambda_1} - E[\bar{\zeta}_{\Lambda_1}|\bar{\eta}_{\Lambda_1}])(\bar{\xi}_{\Lambda_2} - E[\bar{\xi}_{\Lambda_2}|\bar{\eta}_{\Lambda_2}]).$$

Hence,

$$\langle (W_K - F_1 - F_2 + F)^2\rangle = O(K^{-2d}).$$

Now $F_1 + F_2 - 2F$ consist of a term of the form

$$(\bar{\xi}_{\Lambda_1} - E[\bar{\xi}_{\Lambda_1}|\bar{\eta}_{\Lambda_1}])E[\bar{\zeta}_{\Lambda_2}|\bar{\eta}_{\Lambda_2}] - (\bar{\xi}_{\Lambda_2} - E[\bar{\xi}_{\Lambda_2}|\bar{\eta}_{\Lambda_2}])E[\bar{\zeta}_{\Lambda_1}|\bar{\eta}_{\Lambda_1}] \tag{10.19}$$

and another one with the roles of $\xi$ and $\eta$ reversed. Summing by parts, we can rewrite $Av_x\tilde{\Gamma}\tau_x[F_1 + F_2 - 2F]$ as

$$Av_x\phi(\varepsilon x)(\tau_{Ke} - \tau_{-Ke})j(\bar{\xi}_\Lambda - E[\bar{\xi}_\Lambda|\bar{\eta}_\Lambda]), \tag{10.20}$$

where $j = h(\eta_\Lambda) E[\bar{\zeta}_\Lambda | \bar{\eta}_\Lambda]$, plus an analogous term with $\xi$ replaced by $\zeta$. Since we have

$$(10.21) \qquad E[((\tau_{Ke} - \tau_{-Ke})j)^2 (\bar{\xi}_\Lambda - E[\bar{\xi}_\Lambda | \bar{\eta}_\Lambda])^2] = O(K^{-2d}),$$

we have the lemma. $\square$

Now we finish the proof with the replacement of $\hat{F}$ by $F$ (Lemma 10.4). Note here the special role played by the scale $K = O(\varepsilon^{2/d+2})$. We need first a preliminary result on the difference of the chemical potential and the empirical chemical potential.

LEMMA 10.3. *Let $\alpha_x$, $x \in \mathbb{Z}^d$ be a random field taking values in $[-B, B]$ for some $B < \infty$ and satisfying the mixing conditions (2.1). Let $\widehat{\lambda}_K(m)$ be the empirical chemical potential on a block $\Lambda_K$ of side length $K$, given by (8.4). Let $\lambda(m)$ be the annealed chemical potential given by (2.7). Let $\gamma$ be as in (2.1). There is a constant $C$ independent of $m \in (0,1)$ so that*

$$(10.22) \qquad \mathbb{E}[|\widehat{\lambda}_K(m) - \lambda(m)|^\gamma] \leq C K^{-\gamma d/2}.$$

PROOF. Let $I_m$ denote the real interval $[\log \frac{m}{1-m} - B, \log \frac{m}{1-m} + B]$. Since $|\alpha| \leq B$, we have $\lambda(m), \widehat{\lambda}_K(m) \in I_m$. Recall the definition (8.4) of the empirical chemical potential and denote $p(x) = e^x/(1+e^x)$. By Taylor's theorem,

$$m = Av_{x \in \Lambda_K} p(\alpha_x + \lambda) + (\widehat{\lambda}_K - \lambda) Av_{x \in \Lambda_K} p(\alpha_x + \tilde{\lambda})(1 - p(\alpha_x + \tilde{\lambda})),$$

for some $\tilde{\lambda} \in I_m$. Therefore,

$$\mathbb{E}[|\widehat{\lambda}_K - \lambda|^\gamma] \leq r_m^{-1} \mathbb{E}[|Av_{x \in \Lambda_K}(p(\alpha_x + \lambda) - m)|^\gamma],$$

where $r_m = \inf_{y \in I_m} p^\gamma(y)(1-p(y))^\gamma$. By the mixing conditions,

$$\mathbb{E}[|Av_{x \in \Lambda_K}(p(\alpha_x + \lambda) - m)|^\gamma] \leq C K^{-\gamma d/2} \mathbb{E}[|p(\alpha_0 + \lambda) - m|^\gamma].$$

By Taylor's theorem,

$$p(\alpha_0 + \lambda) - m = (\alpha_0 - \alpha^*) p(\tilde{\alpha} + \lambda)(1 - p(\tilde{\alpha} + \lambda)),$$

for some $|\tilde{\alpha}| \leq B$, where $|\alpha^*| \leq B$ is chosen so that $p(\alpha^* - \lambda) = m$. Therefore,

$$\mathbb{E}[|p(\alpha_0 + \lambda) - m|^\gamma] \leq 4 B^\gamma R_m,$$

where $R_m = \sup_{y \in I_m} p^\gamma(y)(1-p(y))^\gamma$. One can check $\sup_{m \in (0,1)} R_m/r_m < \infty$, and, hence, (10.22) follows. $\square$

LEMMA 10.4. *Let*

(10.23) $$\Omega_2 = Av_x\tau_x\Gamma[\hat{F}(m_1,m_2) - F(m_1,m_2)].$$

*Then for any* $\gamma > 0$,

(10.24) $$\limsup_{\delta_2 \to 0}\limsup_{\delta_1 \to 0}\limsup_{\varepsilon \to 0}\sup_f\{\varepsilon^{-1}\langle f\Omega_2\rangle - \gamma\varepsilon^{-2}Av_{x,e}\langle(\nabla_{x,x+e}\sqrt{f}\,)^2\rangle\} \leq 0.$$

PROOF. Write $F - \tilde{F} = K^{-1}[A_1 + A_2]$ with

(10.25) $$A_1 = [(\hat{\lambda}_2 - \lambda)(m_2) - (\hat{\lambda}_1 - \lambda)(m_1)]H,$$

where $H = H(m_1, m_2) = [m_1 e^{-\lambda(m_1)}(1-m_2)e^{\lambda(m_2)} + m_2 e^{-\lambda(m_2)}(1-m_1)e^{\lambda(m_1)}]$
and

(10.26) $$\begin{aligned}A_2 &= \phi(U)b(m_1, m_2) - \phi(-U)b(m_2, m_1),\\ U &= (\hat{\lambda}_2 - \lambda)(m_2) - (\hat{\lambda}_1 - \lambda)(m_1),\end{aligned}$$

where $b(x,y) = xe^{-\lambda(x)}(1-y)e^{\lambda(y)}$ and $\phi(x) = e^x - 1 - x$. Note that $|\hat{\lambda} - \lambda| \leq 2B$ always and, hence, we have $\phi(x), \phi(-x) \leq Cx^2$ for relevant $x$. It is also easy to check that $b$ is uniformly bounded. Hence, we can estimate

(10.27) $$\varepsilon^{-1}K^{-1}Av_x\tau_x\Gamma A_2 \leq C\varepsilon^{-1}K^{-1}Av_x|\hat{\lambda}_K - \lambda|^2,$$

where $\hat{\lambda}_K = \hat{\lambda}_K(\bar{\eta}^K)$ and $\lambda = \lambda(\bar{\eta}^K)$. Next we turn to the $A_1$ term. Summing by parts,

(10.28) $$\varepsilon^{-1}K^{-1}Av_x\tau_x\Gamma H[(\tau_{Ke} - I)(\hat{\lambda}_K - \lambda)] = \varepsilon^{-1}K^{-1}Av_x\tau_x B_x(\hat{\lambda}_K - \lambda).$$

$B = (\tau_{-Ke} - I)\Gamma H$. Let $a > 1$ be as in (2.1) and $a^* = a/(a-1)$ be the conjugate exponent. The right-hand side is bounded above by

(10.29) $$Cq^{a^*}Av_x\tau_x|B|^{a^*} + q^{-a}(\varepsilon K)^{-a}Av_x\tau_x|\hat{\lambda}_K - \lambda|^a.$$

By the two block estimate, for any $C > 0$,

(10.30) $$\limsup_{\varepsilon \to 0}\sup_f\{CAv_x\langle\tau_x|B|^{a^*}f\rangle - \varepsilon^{-2}Av_{x,e}\langle(\nabla_{xx+e}\sqrt{f}\,)^2\rangle\} \leq 0.$$

By Lemma 10.3, and since $K \geq c_0\varepsilon^{-2/d+2}$, there exists $C_0 < \infty$ such that

(10.31) $$\limsup_{\varepsilon \to 0}\mathbb{E}[|(\varepsilon K)^{-1}(\hat{\lambda}_K - \lambda)|^\gamma] \leq C_0.$$

By the mixing conditions,

(10.32) $$\begin{aligned}\mathbb{E}[(Av_x\tau_x((\varepsilon K)^{-1}(\hat{\lambda} - \lambda))^a &- \mathbb{E}[((\varepsilon K)^{-1}(\hat{\lambda} - \lambda))^a])^{\gamma/a}]\\ &\leq C\varepsilon^{\gamma d^2/2a(d+2)},\end{aligned}$$

which is summable in $\varepsilon^{-1} = L = 1, 2, \ldots$ as long as $\gamma > 2(d+2)/d^2$ and $a$ is sufficiently close to 1. By Chebyshev's inequality and the Borel–Cantelli lemma, for almost every realization of the random field,

$$(10.33) \qquad \limsup_{\varepsilon \to 0} q^{-a} Av_x \tau_x((\varepsilon K)^{-1}(\widehat{\lambda} - \lambda)(m^K))^a \leq q^{-a} C_0.$$

Letting $q \to \infty$ completes the proof for $A_1$. The same argument with $\gamma = 2$ and the estimate (10.27) show that $\varepsilon^{-1} K^{-1} Av_x \tau_x \Gamma A_2 \to 0$ with probability one. $\square$

**11. Continuity of the diffusion coefficient.** We introduce a notion of regularity on $D$. There is a finite $C$ such that, for any $x$ and $y$ in $[0,1]$,

$$(11.1) \qquad |D(x) - D(y)|^2 \leq C(x(1-x))^{-1}|x - y|.$$

LEMMA 11.1. *The diffusion coefficient $D(m)$ satisfies* (11.1).

PROOF. Fix a vector $\beta \in \mathbb{R}^d$ and let

$$(11.2) \qquad F(m) = \sqrt{m(1-m)\beta \cdot D^{-1}(m)\beta}.$$

From the proof of the hydrodynamic limit, we have the following representation of the diffusion coefficient:

$$(11.3) \qquad F^2(m) = C \lim_{K \to \infty} K^{-d} \mathbb{E}[\mathbf{V}_K],$$

where

$$(11.4) \quad \mathbf{V}_K = \sup_f \left\{ 2 \sum_e \beta \cdot e \sum_{\substack{x \in \Lambda_K \\ y \in \tau_{Ke} \Lambda_K}} \langle (\eta_y - \eta_x) \nabla_{xy} f \rangle_m - \mathcal{D}_{K \cdot m}(f) \right\},$$

where $\mathcal{D}_{K.m}$ is the Dirichlet form on $\Lambda_K \cup \tau_{Ke} \Lambda_K$. The expectation is with respect to the ergodic invariant measure on that box with density $m$. The constant $C$ comes from the nonstandard average as described in the Introduction. The exact value of $C$ is not relevant.

In order to compare two densities $m$ and $m+h$, we will produce a coupled measure. Independently, at each site $x$, place a red particle with probability $\frac{e^{\alpha_x + \lambda(m)}}{1 + e^{\alpha_x + \lambda(m)}}$, a green particle with probability $\frac{e^{\alpha_x + \lambda(m+h)}}{1 + e^{\alpha_x + \lambda(m+h)}} - \frac{e^{\alpha_x + \lambda(m)}}{1 + e^{\alpha_x + \lambda(m)}}$ and no particle with probability $\frac{1}{1 + e^{\alpha_x + \lambda(m+h)}}$. Let $\rho_x \in \{0, 1\}$ denote the presence or absence of a red particle and $\gamma_x$ the same for a green particle. $\rho_x + \gamma_x = \eta_x$ is the whole configuration which has measure (2.6) with average density $m + h$. We can also couple the dynamics as follows. The red particles have priority in the sense that when a red particle tries to jump on top of green particle, the two particles switch positions. Otherwise, the particles evolve as before. In the common usage, red are *first class particles* and green are *second class particles*. If we are colorblind and see only the total particles, the evolution of $\eta(t)$ is as usual. If we cannot distinguish green particles from empty sites, the evolution of $\rho(t)$ is also as usual. Hence, if $W_K^\rho = K^{-1} \sum_e \beta \cdot e Av_{x \in \Lambda_K, y \in \tau_{Ke}\Lambda_K} w_{xy}^\rho$, where

$$w_{xy}^\rho = \rho_x(1-\rho_y)(1+e^{\alpha_x - \alpha_y}) - \rho_y(1-\rho_x)(1+e^{\alpha_y - \alpha_x})$$

is the red current, we can write

$$(11.5) \qquad F^2(m) = \|W^\rho\|_{m,m+h}^2 = \lim_{K \to \infty} K^{-d} \mathbb{E}[\mathbf{V}_k^\rho],$$

where

$$(11.6) \qquad \mathbf{V}_k^\rho = \sup_f \{2\langle \beta \cdot W_K^\rho f\rangle_{m,m+h} - \mathcal{D}_{K,m,m+h}(f)\},$$

where $\mathcal{D}_{K,m,m+h}$ is the Dirichlet form for the coupled process. Similarly, (11.3) and (11.4) can be rewritten in this language as

$$(11.7) \qquad F^2(m+h) = \|W\|_{m,m+h}^2 = \lim_{K \to \infty} K^{-d} \mathbb{E}[\mathbf{V}_k],$$

where

$$(11.8) \qquad \mathbf{V}_k = \sup_f \{2\langle \beta \cdot W_K f\rangle_{m,m+h} - \mathcal{D}_{K,m,m+h}(f)\}.$$

Let $W^\gamma = W - W^\rho$. We have

$$(11.9) \qquad |F(m+h) - F(m)| = |\|W\| - \|W^\rho\|| \leq \|W^\gamma\|$$

and

$$(11.10) \qquad \|W^\rho\|_{m,m+h}^2 = \lim_{K \to \infty} K^{-d} \mathbb{E}[\mathbf{V}_k^\gamma],$$

where

$$(11.11) \qquad \mathbf{V}_k^\gamma = \sup_f \{2\langle \beta \cdot W_K^\gamma f\rangle_{m,m+h} - \mathcal{D}_{K,m,m+h}(f)\}.$$

Now

$$\langle \beta \cdot W_K^\gamma f\rangle_{m,m+h} = \langle \beta \cdot W_K f\rangle_{m,m+h} - \langle \beta \cdot W_K^\rho f\rangle_{m,m+h}$$
$$= K^{-1} \sum_e \beta \cdot e \; Av_{\substack{x \in \Lambda_K \\ y \in \tau_{Ke}\Lambda_K}} \langle (\gamma_y - \gamma_x)\nabla_{xy}f\rangle_{m,m+h}.$$

By Schwarz's inequality,

$$(11.12) \qquad \begin{aligned} \langle(\gamma_y - \gamma_x)\nabla_{xy}f\rangle &\leq \langle(\gamma_y - \gamma_x)^2\rangle^{1/2}\langle(\nabla_{xy}f)^2\rangle^{1/2} \\ &\leq Ch^{1/2}\langle(\nabla_{xy}f)^2\rangle^{1/2}. \end{aligned}$$

From the moving particles lemma (5.2), we conclude that

$$|F(m+h) - F(m)| \leq C\sqrt{h}. \tag{11.13}$$

From the definition (11.2) of $F$, this gives

$$\sqrt{m(1-m)}|\sqrt{\beta \cdot D^{-1}(m+h)\beta} - \sqrt{\beta \cdot D^{-1}(m)\beta}|$$
$$\leq Ch^{1/2} + |\beta \cdot D^{-1}(m)\beta|^{1/2}|((m+h)(1-(m+h)))^{1/2} - (m(1-m))^{1/2}|.$$

Since $\delta^{-1} I \leq D \leq \delta I$ for some finite $\delta$, this implies that

$$|D(m+h) - D(m)| \leq C(m(1-m))^{-1/2} h^{1/2}. \tag{11.14}$$

By particle–hole duality ($\eta \mapsto 1 - \eta$ and $\alpha \mapsto -\alpha$), we conclude that $|D(m-h) - D(m)| \leq C(m(1-m))^{-1/2} h^{1/2}$ as well. □

Equation (11.1) say that $D(m)$ is Hölder $1/2$ in $(0,1)$, but it says nothing about continuity at the end points. We prove it separately now.

THEOREM 5. *The diffusion coefficient $D(m)$ is continuous on $[0,1]$.*

PROOF. By the symmetry $(\eta, \alpha) \mapsto (1-\eta, -\alpha)$, it follows that $D(m) = \tilde{D}(1-m)$, where $\tilde{D}$ is the diffusion coefficient of the dynamics with $\alpha$ replaced by $-\alpha$. Hence, one only needs to check the continuity at 0.

For any subset $A \in \mathbb{Z}^d$, let

$$\hat{\eta}_A = \prod_{x \in A}(\eta_x - p_x), \qquad \text{where } p_x = \frac{e^{\alpha_x + \lambda(m)}}{1 + e^{\alpha_x + \lambda(m)}}. \tag{11.15}$$

These form an unnormalized orthogonal basis of $L^2(\nu_m)$, where $\nu_m$ is the infinite product measure (2.6) with density $m$. Let $H_n$ denote the subspace spanned by $\hat{\eta}_A$ with $|A| = n$. Let

$$\mathcal{G}_2 = \bigcup_{n \geq 2} H_n. \tag{11.16}$$

We can write the variational formula (2.12) as

$$(\beta, \sigma(m)\beta) = 2 \inf_{U(\alpha)} \left\{ \mathbb{E}\left[\sum_e b_e^2 \langle (\eta_e - \eta_0)^2 \rangle \right] - \Gamma(m) \right\}, \tag{11.17}$$

where

$$b_e = \beta_e + \tau_{-e} U - U, \tag{11.18}$$

$$\Gamma(m) = \sup_{g \in \mathcal{G}_2} \mathbb{E}\left[\sum_e 2b_e \left\langle (\eta_e - \eta_0) \nabla_{0e} \sum_x \tau_x g \right\rangle - \left\langle \left(\nabla_{0e} \sum_x \tau_x g\right)^2 \right\rangle \right]. \tag{11.19}$$

Setting $g = 0$ shows $\Gamma(m) \geq 0$. Hence,

$$(11.20) \quad (\beta, \sigma(m)\beta) \leq 2 \inf_{U(\alpha)} \mathbb{E}\left[\sum_e \langle (\eta_e - \eta_0)^2 \rangle (\beta_e + \tau_{-e} U - U)^2 \right].$$

We now work toward a bound in the other direction. For any $f$,

$$(11.21) \quad \langle (\eta_e - \eta_0) \nabla_{0e} f \rangle = \langle w_{0e} f \rangle$$

and $w_{0,e} \in H_0 \cup H_1 \cup H_2$ with projection onto $H_2$ of the form $c \hat{\eta}_{0,e}$. Let $g = \sum_{|A| \geq 2} \hat{g}_A \hat{\eta}_A$. Note that since $g$ is local, this is a finite sum. Furthermore, note that $\Gamma(m)$ only acts on $\sum_x \tau_x g$, where the sum is over a box large enough that $\nabla_{0,e} \tau_x g = 0$ for any $x$ in its exterior. Because of the sum over shifts, we can assume without loss of generality that $g = \hat{g}_{0,e} \hat{\eta}_{0,e} + f$, where $f \in \mathcal{G}_2$ and $f = \sum_{|A| \geq 2, A \neq \{x, x+e\} \text{ for any } x} \hat{g}_A \hat{\eta}_A$. Then

$$(11.22) \quad \left\langle (\eta_e - \eta_0) \nabla_{0,e} \sum_x \tau_x g \right\rangle = -\hat{g}_{0,e}(p_e - p_0)\langle (\eta_e - \eta_0)^2 \rangle$$

and

$$(11.23) \quad \nabla_{0e} \sum_x \tau_x g = \xi_e + \nabla_{0,e} \sum_x \tau_x f,$$

$$(11.24) \quad \xi_e = (\eta_e - \eta_0)(\hat{\eta}_{-e} \tau_e \hat{g}_{0,e} - (p_e - p_0) \hat{g}_{0,e} + \hat{\eta}_{2e} \tau_{-e} \hat{g}_{0,e}).$$

Hence,

$$(11.25) \quad \Gamma(m) = \sup_{\hat{g}_{0,e}} \left\{ -2 \mathbb{E}\left[ \sum_e b_e \langle (\eta_e - \eta_0)^2 \rangle (p_1 - p_0) \hat{g}_{0e} \right] - \Psi(m) \right\}$$

$$(11.26) \quad \Psi(m) = \inf_f \mathbb{E}\left[ \sum_e \left\langle \left( \xi_e + \nabla_{0,e} \sum_x \tau_x f \right)^2 \right\rangle \right].$$

The supremum is over $\hat{\mathcal{G}}_2 = \{ f \in \mathcal{G}_2 : f \perp \text{span}[\hat{\eta}_{x,x+e}, x \in \mathbb{Z}^d] \}$.

We claim that there is a $C_1 > 0$ such that, for sufficiently small $m > 0$,

$$(11.27) \quad \Psi(m) \geq C_1 m^2 \mathbb{E}[\hat{g}_{0,e}^2].$$

To prove (11.27), note first that

$$(11.28) \quad \xi_e = \xi_e^2 + \xi_e^1 + (\eta_e - \eta_0)(p_e - p_0) \hat{g}_{0,e},$$

where $\xi_e^i \in H_i$. $\hat{\mathcal{G}}_2$ is preserved by shifts, and for $f \in \hat{\mathcal{G}}_2$, $\langle (\eta_e - \eta_0) \nabla_{0,e} \sum_x \tau_x \times f \rangle = 0$. Hence,

$$(11.29) \quad \begin{aligned} \Psi(m) &= \mathbb{E}\left[ \sum_e \langle (\eta_e - \eta_0)^2 \rangle (p_e - p_0)^2 \hat{g}_{0,e}^2 \right] \\ &\quad + \inf_{f \in \hat{\mathcal{G}}_2} \mathbb{E}\left[ \sum_e \left\langle \left( \xi_e^1 + \xi_e^2 - \nabla_{0,e} \sum_x \tau_x f \right)^2 \right\rangle \right]. \end{aligned}$$

Now $\xi_e^2 = (\hat{\eta}_e - \hat{\eta}_0)(\hat{\eta}_{-e}\tau_e \hat{g}_{0,e} + \hat{\eta}_{2e}\tau_{-e}\hat{g}_{0,e})$ and $\xi_e^1 = (p_e - p_0)(\hat{\eta}_{-e}\tau_e \hat{g}_{0,e} + \hat{\eta}_{2e}\tau_{-e}\hat{g}_{0,e})$. One can check then that there is a $c_1 < \infty$ such that

$$(11.30) \quad \mathbb{E}\left[\sum_e \langle (\eta_e - \eta_0)^2 \rangle (p_e - p_0)^2 \hat{g}_{0,e}^2 \right] + \mathbb{E}\left[\sum_e \langle (\xi_e^1)^2 \rangle \right] \leq c_1 m^3 \mathbb{E}[\hat{g}_{0,e}^2].$$

Hence, it is enough to show that, for some fixed $e$, there is a $C_2 > 0$ such that, for small enough $m > 0$,

$$(11.31) \quad \inf_{f \in \hat{G}_2} \mathbb{E}\left[\left\langle \left(\xi_e^2 - \nabla_{0,e} \sum_x \tau_x f \right)^2 \right\rangle \right] \geq C_2 m^2 \mathbb{E}[\hat{g}_{0,e}^2].$$

Without loss of generality, we can write $f = \sum_A \hat{f}_A \hat{\eta}_A$, where $0 \in A$ and $0 \notin \tau_{ne} A$ for $n < 0$. Now one can check using the explicit form of $\xi_e^2$ that $\nabla_{0,e}\tau_x \hat{\eta}_A \perp \xi_e^2$ for all such $A$ except $A \in \{\{0, 2e\}, \{0, e, 2e\}\}$. So the infimum in (11.31) is achieved with $f = \hat{f}_{0,2e}\hat{\eta}_{0,2e} + \hat{f}_{0,e,2e}\hat{\eta}_{0,e,2e}$. Then $\nabla_{0,e}\sum_x \tau_x f$ is given by

$$(11.32) \quad \begin{aligned} & (\hat{f}_{0,2e}\hat{\eta}_{2e} - \tau_{-e}\hat{f}_{0,2e}\hat{\eta}_{3e} - \tau_e \hat{f}_{0,2e}\hat{\eta}_{-e} + \tau_{2e}\hat{f}_{0,2e}\hat{\eta}_{-2e} \\ & \quad - (p_e - p_0)\hat{f}_{0,e,2e}\hat{\eta}_{2e} - (p_e - p_0)\tau_e \hat{f}_{0,e,2e}\hat{\eta}_{-e} \\ & \quad + \tau_{2e}\hat{f}_{0,e,2e}\hat{\eta}_{-2e,-e} + \tau_{-e}\hat{f}_{0,e,2e}\hat{\eta}_{2e,3e}) \\ & \times (\hat{\eta}_e - \hat{\eta}_0 + p_e - p_0) \end{aligned}$$

and we can compute $\mathbb{E}[\langle (\xi_e^2 + \xi_e^1 - \nabla_{0,e}\sum_x \tau_x f)^2 \rangle]$ explicitly to get

$$\begin{aligned} & \mathbb{E}[\{p_e + p_0 - 2p_e p_0\} \\ & \times \{p_{2e}(1 - p_{2e})(\tau_{-e}\hat{g}_{0,e} - \hat{f}_{0,2e} - (p_e - p_0)\hat{f}_{0,e,2e})^2 \\ & \quad + p_{-e}(1 - p_{-e})(\tau_e \hat{g}_{0,e} - \tau_e \hat{f}_{0,2e} - (p_e - p_0)\tau_e \hat{f}_{0,e,2e})^2 \\ & \quad + p_{3e}(1 - p_{3e})(\tau_{-e}\hat{f}_{0,2e})^2 + p_{-2e}(1 - p_{-2e})(\tau_{2e}\hat{f}_{0,2e})^2 \\ & \quad + p_{-2e}(1 - p_{-2e})p_{-e}(1 - p_{-e})(\tau_{2e}\hat{f}_{0,e,2e})^2 \\ & \quad + p_{2e}(1 - p_{2e})p_{3e}(1 - p_{3e})(\tau_{-e}\hat{f}_{0,e,2e})^2 \}]. \end{aligned}$$

Using the bound on the field $\alpha$, one has $p_x \geq C_3 m$ for some $C_3 > 0$. Then it is not hard to check that this quadratic form is bounded below by $C_4 m^2 \mathbb{E}[\hat{g}_{0,e}^2]$ for some $C_4 > 0$. Using the upper bound (11.30) on $\mathbb{E}[\sum_e \langle (\xi_e^1)^2 \rangle]$, we obtain (11.27).

Now (11.27) implies that

$$(11.33) \quad \Gamma(m) \leq C_2 m^2 \mathbb{E}\left[\sum_e b_e^2 \right].$$

Recalling (11.17), we have therefore shown that

$$(11.34) \quad (\beta, \sigma(m)\beta) \geq 2 \inf_{U(\alpha)} \mathbb{E}\left[\sum_e [\langle (\eta_e - \eta_0)^2 \rangle - C_2 m^2](\beta_e + \tau_{-e}U - U)^2\right].$$

Together with the lower bound (11.20) and the formula (2.13) for $D(0)$, we only have to prove that, for any $C$ bounded,

$$(11.35) \quad \lim_{m \to 0} \lambda'(m) 2 \inf_{U(\alpha)} \mathbb{E}\left[\sum_e [\langle (\eta_e - \eta_0)^2 \rangle + C m^2](\beta_e + \tau_{-e}U - U)^2\right]$$
$$= 2z^{-1} \inf_{U(\alpha)} \mathbb{E}\left[\sum_e (e^{\alpha_0} + e^{\alpha_e})(\beta_e + \tau_{-e}U - U)^2\right].$$

Using $\langle (\eta_e - \eta_0)^2 \rangle_m = p_e + p_0 - 2 p_e p_0$ and the form (11.15) of $p_x$, we have

$$(11.36) \quad \langle (\eta_e - \eta_0)^2 \rangle = \frac{e^{\lambda(m)}(e^{\alpha_0} + e^{\alpha_e})}{(1 + e^{\alpha_0 + \lambda(m)})(1 + e^{\alpha_e + \lambda(m)})}.$$

Since $-B \leq \alpha_i \leq B$,

$$(11.37) \quad \frac{e^{\lambda(m)}(e^{\alpha_0} + e^{\alpha_e})}{(1 + e^{B+\lambda(m)})(1 + e^{B+\lambda(m)})} \leq \langle (\eta_e - \eta_0)^2 \rangle$$
$$\leq \frac{e^{\lambda(m)}(e^{\alpha_0} + e^{\alpha_e})}{(1 + e^{-B+\lambda(m)})(1 + e^{-B+\lambda(m)})}.$$

It is elementary to check that as $m \to 0$, $\lambda(m) \sim \log m - \log z$, and $\lambda'(m) \sim m$. Hence, (11.35) follows. $\square$

**12. Uniqueness.** We have the following theorem of uniqueness of weak solutions.

THEOREM 6. *Let $u_1$ and $u_2$ be two weak solutions of $\partial_t u = \nabla \cdot D(u)\nabla u$ [see (2.15)] satisfying (2.14), with the same initial data and suppose that $D$ satisfies (11.1). Then $u_1 = u_2$.*

PROOF. Taking in (2.15), $\phi_\epsilon$ to be the solution at time $\varepsilon > 0$ of the standard heat equation on $\mathbb{T}^d$ with initial data $\delta_0$, and denoting by $f_\varepsilon$ the convolution of a function $f$ with $\phi_\varepsilon$, we have, for $i = 1, 2$,

$$(12.1) \quad \partial_t(u_i)_\varepsilon = \nabla \cdot (D(u_i)\nabla u_i)_\varepsilon$$

as an equality of smooth functions. Let $\psi_\delta(x)$ be an approximation of $|x|$ with $\psi''_\delta(x) = (2\pi\delta)^{-1/2} \exp\{x^2/2\delta\}$. From (12.1), we have

$$\int_{\mathbb{T}^d} \psi_\delta((u_1 - u_2)_\varepsilon(T, \theta)) \, d\theta = -\int_0^T \int_{\mathbb{T}^d} \psi''_\delta((u_1 - u_2)_\varepsilon)(\nabla(u_1 - u_2))_\varepsilon$$
$$\times (D(u_1)\nabla u_1 - D(u_2)\nabla u_2)_\varepsilon \, d\theta \, dt.$$

Let $\varepsilon \to 0$. Since $u_i$ are bounded, $F((u_1 - u_2)_\varepsilon) \to F(u_1 - u_2)$ in $L^\infty$. From (2.14), $\nabla u_i \in L^2 = L^2([0,T] \times \mathbb{T}^d)$ and, since $D$ is bounded, $D(u_i)\nabla u_i \in L^2$ as well. Hence, $(\nabla u_i)_\varepsilon \to \nabla u_i$ and $(D(u_i)\nabla u_i)_\varepsilon \to D(u_i)\nabla u_i$ in $L^2$ and we can take the limit of the above formula to obtain

$$\int_{\mathbb{T}^d} \psi_\delta((u_1 - u_2)(T))\, d\theta$$

$$= -\int_0^T \int_{\mathbb{T}^d} \psi_\delta''(u_1 - u_2)(\nabla u_1 - \nabla u_2)(D(u_1)\nabla u_1 - D(u_2)\nabla u_2)\, d\theta\, dt$$

$$= -\int_0^T \int_{\mathbb{T}^d} \psi_\delta''(u_1 - u_2)(\nabla u_1 - \nabla u_2)D(u_1)(\nabla u_1 - \nabla u_2)\, d\theta\, dt$$

$$+ \int_0^T \int_{\mathbb{T}^d} \psi_\delta''(u_1 - u_2)(\nabla u_1 - \nabla u_2)(D(u_1) - D(u_2))\nabla u_2\, d\theta\, dt.$$

By Schwarz's inequality and the fact that $cI \le D$, we obtain from the previous line that

$$\int_{\mathbb{T}^d} \psi_\delta((u_1 - u_2)(T))\, d\theta \le C \int_0^T \int_{\mathbb{T}^d} \psi_\delta''(u_1 - u_2)|D(u_1) - D(u_2)|^2 |\nabla u_2|^2\, d\theta\, dt$$

for some finite $C$. From (11.1), this is bounded above by

$$(12.2) \qquad C \int_0^T \int_{\mathbb{T}^d} (\psi_\delta''(u_1 - u_2)|u_1 - u_2|)\left(\frac{|\nabla u_2|^2}{u_2(1-u_2)}\right) d\theta\, dt,$$

with perhaps a new finite $C$. Let $F = \frac{|\nabla u_2|^2}{u_2(1-u_2)} \in L^1$ by (2.14). Let $v = u^1 - u^2$ and $h(x) = (2\pi)^{-1/2}|x|\exp\{-x^2/2\}$. Note that $h$ is uniformly bounded with $h(0) = 0$. Rewrite (12.2) as $C \int_0^T \int_{\mathbb{T}^d} h(\delta^{-1/2}v)F\, d\theta\, dt$. The integrand is dominated by $\|h\|_\infty F \in L^1$ and $h(\delta^{-1/2}v) \to 0$ pointwise. Letting $\delta \to 0$ by the dominated convergence theorem, we obtain

$$(12.3) \quad \lim_{\delta \to 0} \int_0^T \int_{\mathbb{T}^d} (\psi_\delta''(u_1 - u_2)|u_1 - u_2|)\left(\frac{|\nabla u_2|^2}{u_2(1-u_2)}\right) d\theta\, dt = 0.$$

One the other hand, by the monotone convergence theorem, we have

$$(12.4) \qquad \lim_{\delta \to 0} \int_{\mathbb{T}^d} \psi_\delta((u_1 - u_2)(T))\, d\theta = \int_{\mathbb{T}^d} |u_1 - u_2|(T)\, d\theta,$$

from which we conclude that $\int_{\mathbb{T}^d} |u_1 - u_2|(T)\, d\theta \le 0$. Hence, $u_1 = u_2$. $\square$

**Acknowledgments.** We would like to thank H. T. Yau for many useful discussions. This work is partly based on notes written jointly with him [18]. We would also like to thank Herbert Spohn for suggesting the problem. The discussion of the structure theorem in Section 7 is partly based on notes of Ana Maria Savu.

DEPARTMENTS OF MATHEMATICS AND STATISTICS
UNIVERSITY OF TORONTO
TORONTO, ONTARIO
CANADA M5S 2E4
E-MAIL: quastel@math.toronto.edu